\documentclass[11pt]{article}

\usepackage[margin=1in]{geometry}
\usepackage{amsmath,amssymb,amsthm,booktabs,array,longtable}
\usepackage[authoryear,round]{natbib}
\usepackage{xcolor}
\usepackage[colorlinks=true,allcolors=blue]{hyperref}
\newtheorem{theorem}{Theorem}[section]
\newtheorem{proposition}[theorem]{Proposition}
\newtheorem{lemma}[theorem]{Lemma}
\newtheorem{corollary}[theorem]{Corollary}

\theoremstyle{definition}
\newtheorem{remark}[theorem]{Remark}
\numberwithin{equation}{section}
\newcommand{\PP}{\mathbb P}
\newcommand{\EE}{\mathbb E}
\newcommand{\RR}{\mathbb R}
\newcommand{\Pbb}{\mathbb P}
\newcommand{\Ebb}{\mathbb E}
\newcommand{\Rbb}{\mathbb R}
\newcommand{\one}{\mathbf 1}
\newcommand{\FDR}{\mathrm{FDR}}
\newcommand{\FDP}{\mathrm{FDP}}

\newcommand{\barPhi}{\overline{\Phi}}

\newcommand{\cC}{\mathcal C}
\newcommand{\cG}{\mathcal G}
\newcommand{\dd}{\,\mathrm d}

\newif\ifstandalonesupplement
\standalonesupplementfalse
\newcommand{\SavedMainSet}[2]{%
  \expandafter\gdef\csname SavedMainNumber@#1\endcsname{#2}}
\SavedMainSet{enh:projection}{2.1}
\SavedMainSet{enh:stability}{3.1}
\SavedMainSet{thm:fixed}{2.2}
\SavedMainSet{thm:nonpositive-main}{4.1}
\SavedMainSet{lem:oracle}{5.1}
\SavedMainSet{thm:one-null}{5.2}
\SavedMainSet{thm:two-null}{5.3}
\SavedMainSet{cor:small-K}{5.4}
\SavedMainSet{prop:eight}{6.1}
\SavedMainSet{cor:K8}{6.2}
\SavedMainSet{thm:truncated}{6.3}
\SavedMainSet{cor:balanced}{6.4}
\SavedMainSet{rankmanual-main}{6.5}
\SavedMainSet{thm:sharp}{6.6}

\newcommand{\MainMap}[3]{%
  \expandafter\def\csname MainLabel@#1\endcsname{#2}%
  \expandafter\def\csname MainKind@#1\endcsname{#3}}
\MainMap{enh:projection}{prop:projection}{Proposition}
\MainMap{enh:stability}{thm:stability}{Theorem}
\MainMap{thm:fixed}{thm:fixed}{Proposition}
\MainMap{thm:nonpositive-main}{thm:nonpositive}{Theorem}
\MainMap{lem:oracle}{lem:oracle}{Lemma}
\MainMap{thm:one-null}{thm:one-null}{Theorem}
\MainMap{thm:two-null}{thm:two-null}{Theorem}
\MainMap{cor:small-K}{cor:small-k}{Corollary}
\MainMap{prop:eight}{prop:eight}{Proposition}
\MainMap{cor:K8}{cor:separation}{Corollary}
\MainMap{thm:truncated}{thm:two-block}{Theorem}
\MainMap{cor:balanced}{cor:balanced}{Corollary}
\MainMap{rankmanual-main}{thm:first-inflation}{Theorem}
\MainMap{thm:sharp}{thm:sharp}{Theorem}
\newcommand{\MainNumber}[1]{%
  \ifstandalonesupplement
    \ifcsname SavedMainNumber@#1\endcsname
      \csname SavedMainNumber@#1\endcsname
    \else ??\fi
  \else
    \ref*{\csname MainLabel@#1\endcsname}%
  \fi}
\DeclareRobustCommand{\MainResult}[1]{%
  \csname MainKind@#1\endcsname~\MainNumber{#1}}
\newwrite\MainNumberFile
\newcommand{\ExportMainOne}[1]{%
  \immediate\write\MainNumberFile{%
    \string\SavedMainSet{#1}{\getrefnumber{\csname MainLabel@#1\endcsname}}}}
\newcommand{\ExportMainNumbers}{%
  \immediate\openout\MainNumberFile=main_labels.tex
  \ExportMainOne{enh:projection}%
  \ExportMainOne{enh:stability}%
  \ExportMainOne{thm:fixed}%
  \ExportMainOne{thm:nonpositive-main}%
  \ExportMainOne{lem:oracle}%
  \ExportMainOne{thm:one-null}%
  \ExportMainOne{thm:two-null}%
  \ExportMainOne{cor:small-K}%
  \ExportMainOne{prop:eight}%
  \ExportMainOne{cor:K8}%
  \ExportMainOne{thm:truncated}%
  \ExportMainOne{cor:balanced}%
  \ExportMainOne{rankmanual-main}%
  \ExportMainOne{thm:sharp}%
  \immediate\closeout\MainNumberFile}

\title{On the false discovery rate of one-sided p-values under Gaussian dependence}
\author{Chenxin Lyu\thanks{Department of Statistics and Actuarial Science, University of Waterloo,  Canada. Email:~{\url{c33lyu@uwaterloo.ca}}}
\and Ruodu Wang\thanks{Department of Statistics and Actuarial Science, University of Waterloo,  Canada. Email:~{\url{wang@uwaterloo.ca}}}}
\date{\today}

\begin{document}
\maketitle

\begin{abstract}
We study the false discovery rate (FDR) of the Benjamini--Hochberg
procedure for $K$ hypotheses at nominal level $q$ under Gaussian
dependence of one-sided p-values.
We ask whether the global null is least favorable after optimizing over
correlation matrices. We show that the global null attains the worst-case
FDR for $K\le3$ at every level. In dimension eight, an explicit mixed-null construction has
FDR $21q/20$ for $0<q<1/2$. It exceeds the global-null maximum
for all sufficiently small $q$, by a uniform
extension of the known fixed-matrix small-level Simes bound.
For $0<q<1/2$, we identify eight as
the first inflation dimension in the monotone rank-one class and
prove that $21q/20$ is its exact eight-dimensional maximum.
In the same level range, the sharp multiplicative bound over all
dimensions in this class is twice the omega constant $w = W(1)$, where
$W$ is the Lambert $W$ function. This bound is approached as the
dimension increases. We also give scale conditions for stability around
the independence benchmark $m_0q/K$, where $m_0$ is the number of true
nulls, and obtain a dimension-free upper bound of $q+q^2/2+O(q^3)$ at
small levels under the global null with nonpositive correlations.
\end{abstract}

\section{Introduction}\label{sec:introduction}

The Benjamini--Hochberg (BH) procedure \citep{BH1995} tests $K$
hypotheses by comparing their ordered p-values with the thresholds
$q/K,2q/K,\ldots,q$. Its error criterion is the false discovery rate
(FDR), the expected proportion of true nulls among the rejected
hypotheses, with this proportion set to zero when there is no
rejection. Under independence and uniform true-null p-values, its
FDR equals $m_0q/K$, where $m_0$ is the number of true nulls.
Positive regression dependence of the full p-vector on the true-null
subset gives the same upper bound \citep{BY2001}. Under arbitrary dependence, the general
bound is $m_0qH_K/K$, where $H_K=\sum_{j=1}^K1/j$
\citep{BY2001,GuoRao2008}.

Gaussian
dependence is more structured than arbitrary dependence, but it
allows both positive and negative correlations.
The purpose of this paper is to study how much
this structure restricts the FDR and which configurations of null versus alternative yield the largest FDR.
Admittedly, most of our results are obtained in low dimensions, making their practical value not very strong because BH is more popular when the number of hypotheses is large. Nevertheless,   they offer theoretical insights on the important problem of FDR control under Gaussian dependence, which is arguably the most standard setting in multiple testing.

We specify dependence through a centered Gaussian vector
$X=(X_1,\ldots,X_K)$ with unit marginal variances and correlation
matrix $\Sigma$, allowing singular matrices. A true-null p-value is
$\Phi(X_i)$, where $\Phi$ is the standard normal distribution
function. At a false null, the p-value is $a_i(X_i)$ for an arbitrary
Borel nondecreasing map $a_i:\mathbb R\to[0,1]$. The matrix and the
maps may depend on $q$. This model separates the Gaussian dependence
from the alternative margins. It includes standard one-sided Gaussian
mean tests with true-null means equal to zero, but also permits
nonparametric and discontinuous alternative margins. Standard
two-sided Gaussian p-values are generally outside this monotone class.

Under the global null, every rejection is false. The FDR then equals
the probability that the Simes function \citep{Simes1986} of the p-values is at most
$q$. Let $s_K(q)$ denote the supremum of this probability over all
$K$-dimensional Gaussian correlation matrices. Let $f_K(q)$ denote
the supremum of FDR when we also optimize over the true-null set and
the alternative maps in the preceding model. Since the global null
is one admissible configuration, $s_K(q)\le f_K(q)$. The question is
whether equality holds. There are two competing effects: rejected
alternatives reduce the proportion of false discoveries, but they
also raise the BH rejection threshold and may cause more true nulls
to be rejected. Thus, allowing alternatives need not reduce the FDR.

We show that the global null is least favorable in dimensions
at most three, but that it need not be least favorable in dimension
eight. Writing $f_{m_0,K}(q)$ for the supremum with exactly $m_0$
true nulls, Theorems~\ref{thm:one-null} and~\ref{thm:two-null} give
\[
 f_{1,K}(q)=\frac{qH_K}{K},\qquad
 f_{2,K}(q)\le s_K(q).
\]
These imply $f_K(q)=s_K(q)$ for $K\le3$ and every $0<q<1$
(Corollary~\ref{cor:small-k}). This identity compares two suprema;
it does not assert that either is bounded by $q$.
In dimension eight, Proposition~\ref{prop:eight} constructs six true
nulls and two monotone alternatives with FDR $21q/20$ throughout
$0<q<1/2$. Combining this construction with the uniform global-null
bound gives
\[
 f_8(q)>s_8(q)\qquad\text{for all sufficiently small }q
\]
(Corollary~\ref{cor:separation}). The mixed configuration therefore
exceeds every global-null Gaussian configuration in the same
dimension. The construction uses general monotone alternatives;
it is not a standard Gaussian mean-shift example.

The small-level Simes upper bound for each fixed Gaussian
correlation matrix is already established by \citet[Theorem~8(III)]{ChiRamdasWang2024}. For the preceding
comparison, we need an error bound that is uniform over these
matrices. Proposition~\ref{thm:fixed} supplies this technical extension:
for every fixed $K$,
\[
 s_K(q)=q\{1+o(1)\},\qquad q\downarrow0.
\]
Here the supremum is taken before the limit, so a nearly maximizing
matrix may change with $q$. The additional argument groups nearly
perfectly positively correlated coordinates and applies Simes
within each group. Correlations between different groups are then
bounded away from one by a constant depending only on $K$, giving
a common remainder bound.

The construction also leads to exact results under rank-one
dependence. Every standardized rank-one Gaussian vector has
coordinates $Z$ or $-Z$, so its true-null p-values form two blocks,
$U=\Phi(Z)$ and $1-U$. The alternatives are monotone in one of the
two directions. This restriction makes it possible to track how
their rejection changes the FDP in each tail. It also has a useful
reduction property: Proposition~\ref{prop:projection} shows that,
for a fixed true-null Gaussian law, the FDR supremum is unchanged
when the alternative latent coordinates are restricted to its
linear span. Our rank-one bounds thus apply whenever the true-null
block has rank one, even if the original alternatives contain
additional Gaussian noise.

Let $f_K^{(1)}(q)$ be the mixed-null supremum restricted to rank-one
latent covariance. Theorems~\ref{thm:first-inflation} and~\ref{thm:sharp}
show that, for every $0<q<1/2$,
\[
 f_K^{(1)}(q)=q\quad(1\le K\le7),\qquad
 f_8^{(1)}(q)=\frac{21}{20}q,\qquad
 \sup_{K\ge1}\frac{f_K^{(1)}(q)}q=2W(1).
\]
Here $W(1)$ is the positive solution of $w e^w=1$, and
$2W(1)=1.1342865808\ldots$. The last constant is approached as
$K\to\infty$, uniformly over $0<q<1/2$, but is not attained in
any finite dimension. The proofs first allow arbitrary alternatives
alongside $a$ copies of $U$ and $b$ copies of $1-U$. For this larger
class, Theorem~\ref{thm:two-block} gives an attained upper bound
expressed through finite harmonic sums, and
Corollary~\ref{cor:balanced} shows that a fixed number of true nulls
should be divided as evenly as possible between the two blocks.
Monotonicity then yields the sharper small-dimensional
classification. The first inflation dimension eight and the exact
value $21q/20$ refer to the monotone rank-one class, not to arbitrary
Gaussian dependence. The unrestricted comparison $f_K=s_K$ remains
open in dimensions four through seven.

Two further bounds describe the effect of dependence.
First, Theorem~\ref{thm:stability} gives an explicit
bound on $|\FDR-m_0q/K|$, uniform over all admissible alternative
maps. For fixed $K\ge2$, write
$d(q)=\min_{i\in I_0,j\ne i}(1-\rho_{ij}(q)^2)$. If
$d(q)\log(1/q)\to\infty$, then $\FDR/q\to m_0/K$ uniformly,
even when the matrix varies with $q$ and approaches perfect absolute
correlation. A fixed positive excess for all sufficiently small $q$
therefore requires at least one pair involving a true null to have
\[
 1-\rho_{ij}^{2}=O\!\left(\frac1{\log(1/q)}\right).
\]
Here \(\rho_{ij}=\Sigma_{ij}\). The pair \((i,j)\) may depend on \(q\). This is a necessary condition, not a sufficient one, and it does
not require the whole matrix to approach rank one. It nevertheless
explains why near-perfect positive and negative correlations
deserve attention in the extremal problem.

Second, under the global null with all off-diagonal correlations
nonpositive, Theorem~\ref{thm:nonpositive} gives a finite-level bound
valid in every dimension. For $0<q\le0.1$, it implies
\[
 \sup_{K\ge1}\frac{s_K^{\mathrm{nonpos}}(q)}q
 \le 1+\sum_{k\ge2}\frac{k^{k-1}(k-1)}{2k!}q^{k-1},
\]
where $s_K^{\mathrm{nonpos}}$ restricts the global-null supremum to
these matrices. At $q=0.05$, the rejection probability is at most
$0.05147$, independently of $K$. This refines the bound of \citet[Theorem~2]{ChiRamdasWang2024} for weakly
negatively dependent p-values, using the stronger Gaussian
assumption. At a first Simes crossing beyond the first threshold,
at least two p-values must lie in the last threshold interval.
A Gaussian density comparison controls these short-interval
events. The resulting bound is $q+q^2/2+O(q^3)$; the coefficient
is not claimed to be sharp. This finite-level conclusion requires
nonpositive correlations and does not bound the unrestricted
mixed-sign supremum.

Our results complement several related approaches. The
FDR-linking theorem of \citet{Su2018} bounds mixed-null FDR
through the joint law of the true-null p-values. Our oracle bound
in Lemma~\ref{lem:oracle} uses this principle while retaining the
finite number of available alternatives. \citet[Theorem~2]{ChenWangZhangWu2026} prove the pointwise small-level
Simes limit for one- and two-sided Gaussian tests when every
off-diagonal correlation is strictly between $-1$ and $1$; this
condition still permits singular matrices. Their Theorem~3 also
treats a configuration with a block of identical p-values and
independent remaining coordinates. For two-sided Gaussian tests under
the global null, \citet[Theorem~2.1]{Zhang2026} establishes exact Simes
validity for sufficiently small levels, uniformly over the correlation
matrix at each fixed dimension.

Recently, \citet{Dobriban2026} shows that BH can exceed its nominal FDR for
two-sided Gaussian tests even when all pairwise correlations are positive.
Allowing the number of tests to grow as $q\downarrow0$, \citet{Lei2026}
obtains lower bounds on the worst-case FDR with leading terms
$q\sqrt{\log(1/q)}/\sqrt{\pi}$ for one-sided tests and
$q\sqrt{\log(1/q)}/(2\sqrt{\pi})$ for two-sided tests.

Our attention is confined to Gaussian dependence, and bounds on FDR and
Simes inequality for arbitrary dependence are largely different; see
\citet{VWW22} and \citet{CLTW23} for combining p-values under arbitrary
dependence and the asymptotic behavior of the Simes inequality in such
settings. Compared with results in our paper and those of \cite{Dobriban2026} and \cite{Lei2026}, Gaussian dependence reduces the upper bound on
FDR substantially compared with arbitrary dependence.

Section~\ref{sec:model} gives the model, the projection reduction,
and the uniform global-null benchmark.
Section~\ref{sec:stability} establishes the stability bound, and
Section~\ref{sec:nonpositive} gives the finite-level bound under
nonpositive correlations. Section~\ref{sec:mixed} treats one and two
true nulls. Section~\ref{sec:rank-one} gives the explicit separation
and sharp rank-one bounds. Short proofs appear with their statements.
The appendices contain the longer proofs and supporting results.

\section{Model and extremal quantities}\label{sec:model}
We describe the dependence of the p-values through a Gaussian latent vector and specify their marginal distributions separately. Under a true null, the p-value is uniform on $[0,1]$. Under a false null, we allow any nondecreasing transformation of the corresponding Gaussian coordinate. This gives a broader class of alternative margins than standard one-sided Gaussian mean tests. Within this framework, we ask whether allowing false nulls can produce a larger FDR than the maximum attainable under the global null in the same dimension.

\subsection{Gaussian model and FDR benchmarks}
Fix $K\ge1$ and $q\in(0,1)$. The BH procedure \citep{BH1995} uses $q$ as its nominal FDR level. For $p\in[0,1]^K$, let
$p_{(1)}\le\cdots\le p_{(K)}$ be its ordered coordinates. Define
\[
 R(p)=\max\{r\in\{1,\ldots,K\}:p_{(r)}\le rq/K\},
\]
with $R(p)=0$ if the set is empty. BH rejects the hypotheses with the
$R(p)$ smallest p-values. If $F$ is the number of rejected true nulls, then
\[
 \FDP=\frac{F}{R\vee1},\qquad \FDR=\EE[\FDP].
\]

Let $\cC_K$ be the set of $K\times K$ correlation matrices, including
singular matrices. Let $X\sim N(0,\Sigma)$, where $\Sigma\in\cC_K$.
Write $\Phi$ and $\phi$ for the standard normal distribution function
and density. For a true-null set $I_0$ with $m_0=|I_0|\ge1$, let
\[
 P_i=
 \begin{cases}
  \Phi(X_i),&i\in I_0,\\
  a_i(X_i),&i\notin I_0,
 \end{cases}
\]
where each $a_i:\RR\to[0,1]$ is Borel measurable and nondecreasing.
The matrix and the alternative maps may depend on $q$.

This framework includes standard one-sided Gaussian mean tests
with true-null means equal to zero. To see this, let
$T\sim N(\mu,\Sigma)$, where $\Sigma_{ii}=1$, and consider the
right-tail p-values $P_i=1-\Phi(T_i)$. Define
$X_i=-(T_i-\mu_i)$. Then $X\sim N(0,\Sigma)$, and the symmetry
of the standard normal distribution gives
\[
  P_i=1-\Phi(T_i)=\Phi(X_i-\mu_i).
\]
Thus, a standard Gaussian mean test corresponds to the
nondecreasing map $a_i(x)=\Phi(x-\mu_i)$. Under a true null,
$\mu_i=0$, so this reduces to $P_i=\Phi(X_i)$. Under an
alternative with $\mu_i>0$, we have
$\Phi(x-\mu_i)<\Phi(x)$ for every $x$, reflecting the tendency
toward smaller p-values in the presence of a positive signal.

More generally, the alternative maps may have jumps or flat regions;
their margins need not belong to a parametric family. Each map moves
in the same direction as $\Phi$ when its latent coordinate increases.
Standard two-sided Gaussian p-values are not generally covered, since
they are not monotone functions of the underlying Gaussian statistic.

The Simes function \citep{Simes1986} and the global-null supremum are
\begin{equation}\label{eq:def-s}
 S_K(p)=\min_{1\le r\le K}\frac{Kp_{(r)}}r,\qquad
 s_K(q)=\sup_{\Sigma\in\cC_K}
 \PP_\Sigma\{S_K(\Phi(X_1),\ldots,\Phi(X_K))\le q\}.
\end{equation}
Let $f_{m_0,K}(q)$ be the supremum of FDR over the preceding model with
exactly $m_0$ true nulls, and put
\[
 f_K(q)=\max_{1\le m_0\le K} f_{m_0,K}(q).
\]
Writing $H_K=\sum_{r=1}^K1/r$, the basic relations are
\[
 q\le s_K(q)\le \min\{H_Kq,1\},\qquad
 f_{K,K}(q)=s_K(q)\le f_K(q).
\]
The lower bound follows from independence, and the upper bound is
the arbitrary-dependence Simes bound \citep{Hommel1983,BY2001}.
The supremum defining $s_K(q)$ is attained. The boundary of the
Simes event is contained in a finite union of marginal threshold
events, each of probability zero. Weak convergence of the Gaussian
law as $\Sigma$ varies therefore gives continuity of the rejection
probability on the compact set $\cC_K$.
Under the global null, every rejection is false, so the FDP equals
$\mathbf 1\{R>0\}$ and the FDR is the Simes rejection probability.
The comparison $f_K(q)=s_K(q)$ therefore asks whether the global null
is least favorable after optimizing the dependence structure.
Let $f^{(1)}_{m_0,K}(q)$ and $f^{(1)}_K(q)$ denote the corresponding
suprema restricted to $\operatorname{rank}(\Sigma)=1$.

\subsection{Projection onto the true-null span}
The true-null law determines the part of the Gaussian dependence that
matters for the extremal comparison. The following reduction removes
residual Gaussian noise from the alternatives without changing the
supremum over their monotone margins.

\begin{proposition}\label{prop:projection}
Fix the Gaussian law of $X_0=(X_i:i\in I_0)$, with covariance rank $r$.
The supremum of FDR over all joint Gaussian extensions and all
nondecreasing alternative maps is unchanged if every alternative latent
coordinate is required to lie in the linear span of $X_0$.
The full latent covariance can therefore be restricted to rank $r$.
\end{proposition}
\begin{proof}
Write $A=\operatorname{Cov}(X_0)$. By Gaussian projection, we can decompose each alternative latent coordinate as $Y_j=b_j^{\mathsf T}X_0+E_j$, where $X_0$ denotes the latent Gaussian vector corresponding to the true nulls, and the residual vector $E$ is independent of $X_0$.
This also holds when $A$ is singular: one may take
$b_j=A^\dagger\operatorname{Cov}(X_0,Y_j)$, since the covariance
vector lies in the range of $A$ and Gaussian uncorrelatedness implies
independence. Here $A^\dagger$ is the Moore--Penrose inverse. Hence, we can rewrite the p-value of a false null as
$P_j=a_j(b_j^{\mathsf T}X_0+E_j)$. Once we condition on $E=e$,
the p-value becomes $a_j(b_j^{\mathsf T}X_0+e_j)$.
Let $\sigma_j^2=\operatorname{Var}(b_j^{\mathsf T}X_0)$.
If $\sigma_j>0$, define
\[
  W_j=\frac{b_j^{\mathsf T}X_0}{\sigma_j}.
\]
Then $W_j$ is standard normal and lies in the linear span
of $X_0$. The conditional p-value is a nondecreasing function
of $W_j$, given by $w\mapsto a_j(\sigma_j w+e_j)$.
If $\sigma_j=0$, the conditional p-value is constant and can
be represented by a constant map of any true-null coordinate.
Thus, the conditional model remains in our admissible class,
with all alternative latent coordinates in the span of $X_0$.

Since $E$ is independent of $X_0$, conditioning does not change
the true-null law. The original FDR is an average of the
conditional FDRs, each of which is bounded by the supremum
over this restricted class. Conversely, the restricted class
is contained in the original class, so the two suprema are equal.
No attainment or measurable choice of an optimizer is needed.
\end{proof}

In particular, if the true-null covariance has rank one,
the optimization reduces to a single Gaussian factor.
The rank-one bounds established below therefore also apply
when only the true-null block has rank one and the alternative
coordinates contain additional Gaussian residual noise.
The reduction concerns a supremum over arbitrary monotone maps;
the conditional rescaling and translation need not preserve the
standard Gaussian mean-shift family.

\subsection{A uniform global-null benchmark}\label{sec:global}
For each fixed Gaussian correlation matrix $\Sigma$, \citet[Theorem~8(III)]{ChiRamdasWang2024} already show that
the global-null Simes rejection probability is at most $q+o(q)$
as $q\downarrow0$. Their result includes singular matrices.
The issue here is uniformity: the remainder for a fixed matrix
does not by itself control the supremum defining $s_K(q)$, since
a nearly maximizing matrix may vary with $q$.

The following result gives a common remainder bound over
$\Sigma\in\mathcal C_K$ for each fixed $K$. We use this technical
extension within the jointly Gaussian class to compare the
eight-dimensional mixed-null construction with the global-null
supremum in Corollary~\ref{cor:separation}.

\begin{proposition}\label{thm:fixed}
For every fixed $K$,
\[
 \lim_{q\downarrow0}\sup_{\Sigma\in\cC_K}
 \frac{\PP_\Sigma\{S_K(P)\le q\}}q=1,
 \qquad P_i=\Phi(X_i).
\]
Equivalently, $s_K(q)=q\{1+o(1)\}$ as $q\downarrow0$.
\end{proposition}
To obtain uniformity, the proof groups the coordinates so that correlations within each
group are nonnegative, whereas correlations between different groups
are bounded away from one by a constant depending only on $K$.
For a group of size $k_g$, the Simes inequality bounds the probability
of a crossing at the original thresholds $rq/K$ by $k_gq/K$.
These bounds sum to $q$. Any remaining rejection requires two
p-values from different groups to be at most $q$. Gaussian tail
bounds show that the total probability of these events is $o(q)$,
uniformly over the correlation matrix. Independence gives the
matching lower bound. The full argument is in Appendix~\ref{supp:proof:thm:fixed}.

This uniform bound provides the comparison needed below. A mixed-null
construction with $\FDR\ge(1+\delta)q$, for a fixed $\delta>0$ and
all sufficiently small $q$, must exceed $s_K(q)$ in that same fixed
dimension once $q$ is small enough. Exceeding $q$ alone would not
establish this comparison.

\begin{remark}\label{rem:fixed}
Proposition~\ref{thm:fixed} is a fixed-dimensional statement; no uniformity in \(K\) is asserted. It does not address the growing-dimensional mixed-null inflation studied by \citet{Lei2026}.
\end{remark}

\section{Stability away from perfect correlation}\label{sec:stability}
Under independence, the BH procedure has FDR equal to $m_0q/K$ when the true-null p-values are uniform \citep{BH1995,BY2001}. This provides a natural benchmark for measuring the effect of dependence, as it accounts for the proportion of true nulls. We therefore ask when Gaussian dependence can produce a fixed positive excess of $\FDR/q$ over $m_0/K$ as $q\downarrow0$ with $K$ fixed.

The relevant quantity is the smallest departure from perfect absolute
correlation among pairs involving a true null:
\[
 d(\Sigma,I_0)=\min_{i\in I_0,\,j\ne i}(1-\rho_{ij}^2),
 \qquad \rho_{ij}=\Sigma_{ij}.
\]
The next result allows this quantity to tend to zero. It is enough
that it decay more slowly than $1/\log(1/q)$; neither the correlation
matrix nor the alternative maps need to remain fixed.

\begin{theorem}\label{thm:stability}
Fix $K\ge2$ and a nonempty true-null set $I_0$, with $m_0=|I_0|$.
Let $\eta(q)\in(0,1]$ satisfy
$\eta(q)\log(1/q)\to\infty$ as $q\downarrow0$. Then
\[
 \sup_{\substack{\Sigma\in\cC_K:\ d(\Sigma,I_0)\ge\eta(q)\\
                  a_j:\RR\to[0,1]\ \text{Borel nondecreasing},\ j\notin I_0}}
 \left|\frac{\FDR}{q}-\frac{m_0}{K}\right|\longrightarrow0.
\]
For any admissible family of correlation matrices $\Sigma(q)$ and
alternative maps, without imposing the preceding lower bound on $d$,
if $\FDR/q\ge m_0/K+\delta$ for some fixed $\delta>0$
and all sufficiently small $q$, then
\[
 d(\Sigma(q),I_0)=O\!\left(\frac1{\log(1/q)}\right).
\]
\end{theorem}

The supremum allows both the matrix and the alternative maps to
depend on $q$. For example, $d(q)=\{\log(1/q)\}^{-1/2}$ is allowed,
even though the largest absolute correlation involving a true null
then tends to one. The condition concerns only pairs involving a
true null; it does not require the full matrix to be nonsingular.

The second conclusion says that a persistent positive excess requires
at least one pair involving a true null to approach perfect absolute
correlation at the stated scale. The pair may vary with $q$.
This is a necessary condition, not a sufficient one, and it does not
require the whole matrix to approach rank one.

Both conclusions follow from the following finite-level bound. Whenever $d=d(\Sigma,I_0)>0$ and $0<q<1/2$,
\[
 \left|\FDR-\frac{m_0q}{K}\right|
 \le \frac{(K-1)M(q)}{\sqrt{2\pi}}
 \sum_{i\in I_0}\sum_{j\ne i}
 \frac{|\rho_{ij}|}{\sqrt{1-\rho_{ij}^2}},
 \qquad M(q)=\phi(t)-tq,\quad t=\barPhi^{-1}(q),
\]
where $\barPhi=1-\Phi$. In particular,
\[
 \left|\frac{\FDR}{q}-\frac{m_0}{K}\right|
 \le \frac{C^{\mathrm{mon}}_{K,m_0}}{t\sqrt d},
 \qquad C^{\mathrm{mon}}_{K,m_0}=\frac{m_0(K-1)^2}{\sqrt{2\pi}}.
\]

Since $t^2\sim2\log(1/q)$, the last bound proves the uniform convergence.
A fixed excess $\FDR/q\ge m_0/K+\delta$ also requires
\[
 d\le \min\left\{1,\frac{(C^{\mathrm{mon}}_{K,m_0})^2}{\delta^2t^2}\right\},
\]
which gives the necessary scale. This condition is automatic when $d=0$.

The proof compares the FDR contribution of each true null with $q/K$.
Temporarily replacing that null's p-value by zero gives a rejection
count that expresses its contribution. Writing $T=-X$, Gaussian
regression separates the remaining coordinates into linear functions
of $T_i$ and a residual vector independent of $T_i$. Freezing the
linear terms at $t$ gives a reference count independent of the
null p-value, whose contribution is exactly $q/K$. The error requires
one of the BH threshold comparisons to change. The relevant residual
has variance $1-\rho_{ij}^2$, so its density bound deteriorates near
perfect absolute correlation. Integrating the threshold-crossing
bound over $T_i\ge t$ produces $M(q)$.

The proof is in Appendix~\ref{supp:proof:thm:stability}.

\section{Bounds under nonpositive correlations}
\label{sec:nonpositive}\label{sec:sign-sectors}

Proposition~\ref{thm:fixed} describes the limit as $q$ tends to zero.
To obtain bounds at specified testing levels, we now restrict the
signs of the correlations. Under the global null, nonnegative
Gaussian correlations imply positive regression dependence, which
gives $\PP_\Sigma\{S_K(P)\le q\}\le q$ \citep{BY2001}.
With nonpositive off-diagonal correlations, the Simes inequality
reverses \citep{BlockSavitsWang2008}. Strict inflation can occur
already for two negatively correlated one-sided Gaussian tests
\citep{SamuelCahn1996}.

\begin{theorem}\label{thm:nonpositive}\label{thm:nonpositive-main}
Under the global null, suppose $\Sigma_{ij}\le0$ for $i\ne j$ and
$0<q\le0.1$. Then
\begin{equation}\label{eq:nonpositive-main}
 \PP_\Sigma\{S_K(P)\le q\}
 \le q+\sum_{k=2}^K\binom Kk\binom k2
       \left(\frac qK\right)^2
       \left(\frac{kq}{K}\right)^{k-2}.
\end{equation}
The conclusion includes singular correlation matrices.
If $s_K^{\mathrm{nonpos}}$ denotes the supremum in \eqref{eq:def-s}
restricted to these matrices, then
\begin{equation}\label{eq:nonpositive-dimensionfree}
 \sup_{K\ge1}\frac{s_K^{\mathrm{nonpos}}(q)}q
 \le 1+\sum_{k\ge2}\frac{k^{k-1}(k-1)}{2k!}q^{k-1}.
\end{equation}
\end{theorem}

\begin{proof}
First suppose $\Sigma$ is positive definite. We keep track of the first
Simes crossing.
Put $h=q/K$, define $N(u)=\sum_{i=1}^K\mathbf1\{P_i\le u\}$, and let
\[
 T=\min\{k\in\{1,\ldots,K\}:N(kh)\ge k\},
\]
with $T=\infty$ if there is no crossing. BH rejects whenever $T\le K$,
although its rejection count is the largest crossing index.
The event $T=1$ has probability at most $Kh=q$. If $T=k\ge2$, then
\[
 N((k-1)h)\le k-2,\qquad N(kh)\ge k.
\]
Hence at least two p-values lie in the last interval
$I_k=((k-1)h,kh]$. For instance, a first crossing at $k=3$ needs
three p-values below $3h$, but at most one can be below $2h$;
at least two must lie in $(2h,3h]$.

The Gaussian assumption makes this interval restriction useful.
By Lemma~\ref{supp:lem:density-dom},
the joint density of any selected set of p-values is at most one
when all of them lie in $(0,q]$. Thus, in this region, the probability
of a rectangle is bounded by its volume under independent uniform
coordinates. The Gaussian density comparison holds for latent
coordinates at most $-1$, which covers the required region because
$q\le0.1<\Phi(-1)$.

For a fixed set of $k$ coordinates and a selected pair, the event that
the pair lies in $I_k$ and the remaining $k-2$ coordinates lie in
$[0,kh]$ therefore has probability at most $h^2(kh)^{k-2}$.
There are $\binom Kk$ choices of the set and $\binom k2$ choices
of the pair. A union bound, followed by summation over $T$, gives
\eqref{eq:nonpositive-main}.

For singular $\Sigma$, apply the bound to
$\Sigma_\delta=(1-\delta)\Sigma+\delta I$. The signs are preserved,
and the Gaussian laws converge weakly as $\delta\downarrow0$.
The Simes-event boundary lies in a finite union of marginal grid levels,
each of probability zero, so the bound passes to the limit.

Finally, divide \eqref{eq:nonpositive-main} by $q$, use
$\binom Kk/K^k\le1/k!$, and extend the nonnegative sum to all $k\ge2$
to obtain \eqref{eq:nonpositive-dimensionfree}.
\end{proof}

At $q=0.01,0.05,0.10$, the right-hand side of
\eqref{eq:nonpositive-dimensionfree} is approximately $1.005154$,
$1.029325$, and $1.070406$, respectively. For example,
at the nominal level $q=0.05$, the rejection probability is at most
$0.05147$ for every $K$. These are upper bounds, not claimed maxima.

For comparison, \citet[Theorem~2]{ChiRamdasWang2024} prove that weakly negatively
dependent uniform p-values satisfy
\[
 \PP\{S_K(P)\le q\}
 \le q+\sum_{k=2}^K\binom Kk\left(\frac{kq}{K}\right)^k,
 \qquad 0<q<1.
\]
Their assumption requires
$\PP(P_i\le u\text{ for all }i\in A)\le u^{|A|}$ for every
nonempty index set $A$ and $u\in(0,1)$. It includes the nonpositive
Gaussian class but does not require Gaussianity. Our additional
density comparison controls short-interval events, which do not
follow from those lower-orthant bounds alone.

For each $k\ge2$, the ratio of the corresponding summands in the
two bounds is
\[
 \frac{\binom k2 h^2(kh)^{k-2}}{(kh)^k}
 =\frac{\binom k2}{k^2}=\frac{k-1}{2k}.
\]
This coefficient has a simple counting interpretation. Consider $k$
independent uniform points on $[0,kh]$ only as a reference for
measuring volume. The number $B$ in the last interval $I_k$ has
distribution $\operatorname{Binomial}(k,1/k)$, and
\[
 \PP(B\ge2)\le\EE\binom B2
 =\frac{k-1}{2k}.
\]
Whenever at least two points occupy the last interval, there is at
least one unordered pair there. Counting pairs may overcount
configurations with three or more points, which explains the
inequality. The reference calculation does not assert independence
or a binomial occupancy law for the actual p-values.

The coefficient equals $1/4$ at $k=2$ and approaches $1/2$ as $k$
increases. In particular, the leading correction in the dimension-free
probability bound decreases from $2q^2$ to $q^2/2$, while the order
$q+O(q^2)$ is unchanged. This compares upper bounds and does not
claim sharp coefficients or a halving of the actual rejection
probability. The result concerns the global null with nonpositive
correlations; arbitrary mixed-sign correlations and mixed-null FDR
are outside its scope.

\section{Mixed-null reductions}\label{sec:mixed}

We now compare mixed-null FDR with the global-null supremum.
We observe that rejecting false nulls lowers the proportion of false discoveries, but
it can also increase the BH threshold and cause additional true nulls
to be rejected. To control these two effects, we first fix the
true-null p-values and optimize over the number of rejected false nulls.
This gives a bound that depends only on the true-null law.

\subsection{An oracle bound}

Let $n=m_0$, let $\ell=K-n$ be the number of false nulls, and write
$Q_{(1)}\le\cdots\le Q_{(n)}$ for the ordered true-null p-values,
with $Q_{(0)}=0$. For each hypothetical number
$m\in\{0,\ldots,\ell\}$ of rejected false nulls, define
\[
 j_m(Q)=\max\left\{0\le j\le n:
 Q_{(j)}\le\frac{q(m+j)}K\right\}.
\]
Thus, $j_m(Q)$ is the largest number of true nulls compatible with
the BH threshold when $m$ false nulls are rejected. Set
\begin{equation}\label{eq:oracle-H}
 H_{n,K}(Q)=\max_{0\le m\le\ell}
 \frac{j_m(Q)}{m+j_m(Q)},
 \qquad 0/0=0 \text{ by convention}.
\end{equation}
The restriction $m\le K-n$ retains the number of alternatives
available in the original problem. We use the FDR-linking approach of \citet[Section~2.1]{Su2018}, while keeping this finite constraint
on the number of rejected false nulls.

\begin{lemma}\label{lem:oracle}
For every realized full p-vector,
\[
 \FDP\le H_{n,K}(Q).
\]
Consequently,
\[
 \FDR\le\EE H_{n,K}(Q).
\]
\end{lemma}

\begin{proof}
If BH rejects $m$ false nulls and $j\ge1$ true nulls, then
$Q_{(j)}\le q(m+j)/K$, so $j\le j_m(Q)$. Since
$x/(m+x)$ is nondecreasing in $x$,
\[
 \FDP=\frac{j}{m+j}
 \le\frac{j_m(Q)}{m+j_m(Q)}\le H_{n,K}(Q).
\]
The case $j=0$ is immediate. Taking expectations proves the FDR bound.
\end{proof}

Taking the supremum over Gaussian true-null laws gives the corresponding
bound on $f_{n,K}(q)$. The maximizing choice in
\eqref{eq:oracle-H} may depend on the entire realized vector $Q$.
Consequently, this pointwise upper bound need not be attained by
admissible nondecreasing Gaussian alternative maps.

\subsection{One and two true nulls}

With one true null, the smallest rejection count that can include it
determines its largest possible contribution to FDR. The resulting
bound is attained by a monotone step construction.
For $K=2$, \citet[Section~4.1.2]{RoquainVillers2011} attain the same
value $3q/4$ for $0<q<1/2$ in a Gaussian mean-shift model.

\begin{theorem}\label{thm:one-null}
For every $K\ge1$ and $0<q<1$,
\[
 f_{1,K}(q)=\frac qK H_K.
\]
\end{theorem}

\begin{proof}
Write $h=q/K$ and let $U$ be the true-null p-value. If
$U\in((r-1)h,rh]$ is rejected, then BH must reject at least $r$
hypotheses, so the FDP is at most $1/r$. If $U>q$, it cannot be
rejected. Uniformity of $U$ therefore gives
\[
 \FDR\le\sum_{r=1}^K\frac{h}{r}=\frac qK H_K.
\]
For equality, take $U=\Phi(Z)$ with $Z\sim N(0,1)$ and define
the $K-1$ false-null p-values by
\[
 A_j(U)=
 \begin{cases}
  0,&U>jh,\\
  1,&U\le jh,
 \end{cases}
 \qquad j=1,\ldots,K-1.
\]
Each $A_j$ is a nondecreasing function of the Gaussian coordinate
$-Z$, so this is an admissible rank-one model. On
$U\in((r-1)h,rh]$, exactly $r-1$ false-null p-values are zero
and the others are one. BH therefore rejects exactly $r$
hypotheses, including the true null, and the FDP equals $1/r$.
The construction supplies exactly enough rejections to reach the
next BH threshold. Integration over the $K$ intervals gives equality;
their endpoints have probability zero.
\end{proof}

For two true nulls, we compare the oracle bound with global-null
vectors formed by repeating the two observations. Each repetition
pattern preserves Gaussian dependence, and averaging over the
patterns gives the required comparison.

\begin{theorem}\label{thm:two-null}
For every $K\ge2$ and $0<q<1$,
\[
 f_{2,K}(q)\le s_K(q).
\]
\end{theorem}

\begin{proof}
Let $U,V$ be the true-null p-values, and put
$x=\min(U,V)$ and $y=\max(U,V)$. For $c=1,\ldots,K-1$, let
$P^{(c)}$ contain $c$ copies of $U$ and $K-c$ copies of $V$.
Every $P^{(c)}$ is an admissible global-null Gaussian vector,
possibly with singular covariance. We claim that
\begin{equation}\label{eq:clone-average}
 H_{2,K}(U,V)\le\frac1{K-1}\sum_{c=1}^{K-1}
 \one\{S_K(P^{(c)})\le q\}.
\end{equation}
If $y\le q$, every indicator is one. If $y>q$, at most one
true null can be rejected. Let
\[
 r=\min\left\{k\in\{1,\ldots,K-1\}:x\le\frac{kq}K\right\}.
\]
When this set is empty, both sides of \eqref{eq:clone-average}
are zero. Otherwise, the oracle uses $r-1$ rejected false nulls
and gives $H_{2,K}(U,V)=1/r$. Across the $K-1$ repetition
patterns, the multiplicity of $x$ runs through $1,\ldots,K-1$.
The Simes event occurs exactly when this multiplicity is at least
$r$. Hence the average in \eqref{eq:clone-average} equals
$(K-r)/(K-1)$, and
\[
 \frac{K-r}{K-1}\ge\frac1r,
 \qquad 1\le r\le K-1,
\]
because $r(K-r)-(K-1)=(r-1)(K-r-1)\ge0$.
Taking expectations in \eqref{eq:clone-average}, each probability
on the right is at most $s_K(q)$. Lemma~\ref{lem:oracle} and a supremum
over the bivariate Gaussian true-null law prove the result.
\end{proof}

These two results exhaust the mixed configurations in dimensions
at most three.

\begin{corollary}\label{cor:small-k}
For $K=1,2,3$ and every $0<q<1$,
\[
 f_K(q)=s_K(q).
\]
\end{corollary}

\begin{proof}
The global-null configuration gives $f_K(q)\ge s_K(q)$.
For $K\ge2$, Theorem~\ref{thm:one-null} gives
$f_{1,K}(q)=qH_K/K<q\le s_K(q)$, and
Theorem~\ref{thm:two-null} controls the only other mixed
configuration when $K\le3$. For $K=1$, the two quantities coincide
by definition.
\end{proof}

In any dimension, a strict separation $f_K(q)>s_K(q)$ must involve between three
and $K-1$ true nulls. The unrestricted comparison remains open
for $K=4,\ldots,7$. The next section gives a strict separation
in dimension eight for sufficiently small $q$.

\section{Sharp bounds in the rank-one model}\label{sec:rank-one}

We now use rank-one dependence to construct a mixed-null configuration
that exceeds the global-null supremum at sufficiently small levels,
and then obtain matching extremal bounds. The construction yields $21q/20$ in dimension eight.
We ask whether this value is optimal, whether inflation is possible in
smaller dimensions, and how large $\FDR/q$ can become as the dimension
increases. A bound with arbitrary alternatives will provide a useful
comparison before we impose monotonicity.

The two-block form follows from the rank-one assumption. Since each
latent coordinate has variance one, a centered rank-one Gaussian vector
can be represented as
\[
 X_i=s_iZ,\qquad s_i\in\{-1,1\},\qquad Z\sim N(0,1).
\]
Put $U=\Phi(Z)$. Each true-null p-value is then either
$\Phi(Z)=U$ or $\Phi(-Z)=1-U$. An alternative with sign $s_j$ is
$a_j(s_j\Phi^{-1}(U))$, which is nondecreasing in $U$ when $s_j=1$
and nonincreasing when $s_j=-1$. Thus the two null blocks and these
two possible directions describe every model in the rank-one class.

\subsection{An explicit separation in dimension eight}\label{sec:eight}

The reductions in Section~\ref{sec:mixed} leave open whether a mixed-null
configuration can have larger FDR than every global-null configuration in
the same dimension. We construct such an example from one Gaussian
variable. A single tied null block cannot produce inflation above $q$:
if the $m$ true-null p-values all equal one uniform variable, then
\[
 \FDR\le\frac{mq}{K}\{1+H_K-H_m\}\le q,
\]
even with arbitrary alternatives (Appendix~\ref{supp:proof:prop:staircase}).
We therefore use two blocks,
equal to $U$ and $1-U$, so that false discoveries can occur in either
tail. The alternatives must increase the BH threshold without reducing
the FDP too much. To achieve this, they become jointly rejectable only
when a null block can no longer be rejected on its own.

Let $U=\Phi(Z)$ for $Z\sim N(0,1)$, and let $0<q<1/2$,
$h=q/8$, and $\lambda=9h/2$. Take
\[
 P=(U,U,U,1-U,1-U,1-U,A_1(U),A_2(U)),
\]
with the first six coordinates true nulls and the last two alternatives, where
\[
 A_1(u)=
 \begin{cases}1,&u\le3h,\\ \lambda,&u>3h,\end{cases}
 \qquad
 A_2(u)=
 \begin{cases}\lambda,&u<1-3h,\\1,&u\ge1-3h.\end{cases}
\]
In the left tail, $A_2$ already equals $\lambda$, but
$\lambda>4h$ prevents it from being rejected together with just the
three nulls. Once $U>3h$, $A_1$ drops to the same value. Since
$\lambda<5h$, the pair makes rank five available and extends null
rejection from $U\le3h$ to $U\le5h$. The right tail is symmetric.

\begin{proposition}\label{prop:eight}
The preceding vector belongs to the monotone rank-one Gaussian class and has
\[
 \FDR=2\left(3h+2h\frac35\right)=\frac{21}{20}q.
\]
\end{proposition}

\begin{proof}
To check admissibility, use $Z$ for the first three latent coordinates
and $-Z$ for the next three. The map $A_1(\Phi(Z))$ is nondecreasing
in $-Z$, and $A_2(\Phi(Z))$ is nondecreasing in $Z$. Thus every
coordinate is an allowed nondecreasing map of its Gaussian coordinate,
and the latent covariance has rank one.

The complete BH calculation is as follows; interval endpoints have
probability zero.
\begin{center}
\begin{tabular}{lc cc}
\toprule
Range of $U$ & Rejected p-values & $R$ & $\FDP$\\
\midrule
$(0,3h)$ & $U,U,U$ & $3$ & $1$\\
$(3h,5h)$ & $U,U,U,A_1,A_2$ & $5$ & $3/5$\\
$(5h,1-5h)$ & None & $0$ & $0$\\
$(1-5h,1-3h)$ & $1-U,1-U,1-U,A_1,A_2$ & $5$ & $3/5$\\
$(1-3h,1)$ & $1-U,1-U,1-U$ & $3$ & $1$\\
\bottomrule
\end{tabular}
\end{center}
Both null blocks cannot be rejected because every rejected p-value is
at most $q<1/2$. In the middle interval, neither block can be rejected
at rank at most five, and the two alternatives alone cannot be rejected
because $\lambda>2h$. Integrating the last column gives the stated
FDR: each tail contributes $3h$ with FDP one and $2h$ with FDP $3/5$.
\end{proof}

\begin{corollary}\label{cor:separation}
There exists $q_0\in(0,1/2)$ such that
\[
 f_8(q)>s_8(q),\qquad 0<q\le q_0.
\]
\end{corollary}

\begin{proof}
Proposition~\ref{prop:eight} gives $f_8(q)\ge21q/20$ for every
$0<q<1/2$. Proposition~\ref{thm:fixed} gives $s_8(q)/q\to1$ as
$q\downarrow0$, so $s_8(q)<21q/20$ for all sufficiently small $q$.
\end{proof}

The example therefore exceeds the nominal level throughout $(0,1/2)$,
whereas its separation from the global-null supremum is established
only at sufficiently small levels.

\subsection{A bound for two countermonotone blocks}
Let $K=a+b+\ell$, where $a,b,\ell$ are nonnegative integers.
The true-null p-values are $a$ copies of $U$ and $b$ copies of $1-U$.
In this subsection the $\ell$ alternative p-values may be arbitrary random
variables in $[0,1]$, with arbitrary dependence on $U$ and each other.
For $m\ge0$, set
\[
 G_\ell(m)=m+\sum_{j=1}^{\ell}\frac{m}{m+j},\qquad G_\ell(0)=0.
\]
\begin{theorem}\label{thm:two-block}
For $0<q<1/2$, every model in this two-block class satisfies
\[
 \frac{\FDR}{q}\le\frac{G_\ell(a)+G_\ell(b)}K.
\]
For each fixed split $(a,b,\ell)$ with $K\ge1$, the bound is attained in this class.
\end{theorem}

\begin{proof}
Put $h=q/K$. The two null blocks cannot both be rejected, since
every rejected p-value is at most $q$ and $U+(1-U)=1>2q$.
If the $U$ block is rejected and $a>0$, all its tied p-values are
rejected and $a\le R\le a+\ell$. Its FDP contribution is therefore at most
\[
 \frac{a}{\max\{a,\lceil U/h\rceil\}}
 \one\{U\le(a+\ell)h\}.
\]
Integration gives $hG_\ell(a)$; the other tail gives $hG_\ell(b)$.
A missing block contributes zero. The bound is pointwise in the
alternatives, so their additional randomness or dependence causes no issue.

For equality, take the $j$th alternative to be zero on
\[
 (a+j-1)h<U<1-(b+j-1)h
\]
and one elsewhere, for $j=1,\ldots,\ell$.
On the relevant left tail $U\le(a+\ell)h\le q$, every upper endpoint
exceeds $1-q$. Thus, if $a>0$, BH rejects only the $a$ nulls on
$0<U\le ah$, and exactly those nulls and $j$ zero alternatives on
$(a+j-1)h<U\le(a+j)h$. The opposite block exceeds $q$ and the
remaining alternatives equal one. This attains each layer of the bound;
the right tail is symmetric. Outside these tails the rank restriction
excludes false discoveries. Missing blocks and $\ell=0$ cause no change,
and interval endpoints have probability zero.
\end{proof}

The cutoffs $a+\ell=K-b$ and $b+\ell=K-a$ matter: the opposite null
block cannot help raise the rejection rank. The attaining alternatives
are generally nonmonotone, so this is not an exact formula for the
original monotone class.

Theorem~\ref{thm:two-block} fixes the sizes of both null blocks.
If only their total size is fixed, we can also optimize its allocation
between the two blocks.

\begin{corollary}\label{cor:balanced}
Let $K\ge1$, $n\in\{1,\ldots,K\}$, and $0<q<1/2$.
In the unrestricted-alternative two-block class with $n$ true nulls,
the maximum FDR is
\[
 \frac qK\left[
 G_{K-n}\bigl(\lfloor n/2\rfloor\bigr)+
 G_{K-n}\bigl(\lceil n/2\rceil\bigr)\right].
\]
The maximum over the true-null count is the maximum of these $K$ quantities.
\end{corollary}

\begin{proof}
Put $\ell=K-n$. Extending $G_\ell$ to real $t\ge0$ gives
\[
 G_\ell''(t)=-2\sum_{j=1}^{\ell}\frac{j}{(t+j)^3}\le0
 \qquad(t\ge0).
\]
For fixed integer $a+b=n$, concavity makes $G_\ell(a)+G_\ell(b)$
largest at the most balanced integer split.
Theorem~\ref{thm:two-block} gives attainability. If $\ell=0$,
every split has the same value, so the formula still holds.
\end{proof}

The case $n=0$ has FDR zero. The maximum in
Corollary~\ref{cor:balanced} allows nonmonotone alternatives.
For example, $(a,b,\ell)=(1,2,1)$ gives $25q/24$ in dimension
four. The next theorem shows why this larger-class value cannot be
attained by monotone Gaussian alternatives.
\subsection{The first inflation dimension}
We return to the original rank-one Gaussian model. Each alternative is
a nondecreasing or nonincreasing function of $U$, so its behavior in the
two tails cannot be chosen separately. In particular, a nonincreasing
alternative that becomes small in the left tail remains small in the
right tail. It may then increase the denominator of the FDP where an
unrestricted construction would prefer to leave it unrejected.

\begin{theorem}\label{thm:first-inflation}
For every $0<q<1/2$,
\[
 f_K^{(1)}(q)=q\quad(1\le K\le7),\qquad
 f_8^{(1)}(q)=\frac{21}{20}q.
\]
\end{theorem}

Write \(n=a+b=m_0\), \(\ell=K-n\), and \(h=q/K\). The upper bound uses this constraint in a small number of cases. In one tail, let $m\ge1$ be the size of its null block and write its
common p-value as $hx$, where $0<hx\le q$. Throughout this tail
the opposite block exceeds $q$ and cannot be rejected. Set the $m$
p-values in the first block to zero, and let
$r(x)$ be the resulting BH rejection count. Restoring them gives the
exact null contribution
\[
 \frac{m}{r(x)}\one\{x\le r(x)\},\qquad
 m\le r(x)\le m+\ell.
\]
If all alternatives are nondecreasing in the tail distance $x$, then
$r(x)$ is nonincreasing. Its null-rejection region is an initial
interval, and integration bounds that tail's FDR by $mh$.
If all but at most one alternative are nondecreasing, the same argument,
applied with the remaining alternative fixed at zero or one, gives
\[
 \FDR_{\rm tail}/h\le P_\ell(m):=
 m+\frac{m\ell}{m+\ell}.
\]
The endpoint arguments allowing discontinuous maps are in
Lemmas~\ref{supp:rankmanual-rank} and~\ref{supp:rankmanual-onehelper}.

A single null block is bounded by $hG_\ell(m)\le q$.
If all alternatives have the same direction, then, after exchanging
the two tails if necessary, one contributes at most $ah$ and the
other at most $hG_\ell(b)$, whose sum is at most
$q$. This includes every case with one alternative. With two
oppositely directed alternatives, each tail satisfies the $P_2$ bound;
with three alternatives of mixed directions, one tail satisfies the
$P_3$ bound and the other the $G_3$ bound. Among splits with two nonempty
null blocks and alternatives of both directions, the unrestricted
bound and balanced-split comparison suffice through dimension eight
except in the following five cases, for which these refinements give
\[
\begin{array}{c|ccccc}
 (K,n,\ell)&(6,4,2)&(7,4,3)&(7,5,2)&(8,5,3)&(8,6,2)\\ \hline
 \text{bound on }\FDR/q&1&29/30&33/35&121/120&21/20.
\end{array}
\]
For example, the last entry uses the bound for two oppositely directed
alternatives. With $a+b=6$, concavity of $P_2$ makes
$P_2(a)+P_2(b)$ largest at $a=b=3$, giving
\[
 \frac{P_2(3)+P_2(3)}8
 =\frac{2(3+6/5)}8=\frac{21}{20}.
\]
For $(K,n,\ell)=(7,5,2)$, the bound also uses the fact that an
alternative made small in one tail remains small in the other;
Lemma~\ref{supp:rankmanual-asymmetric} gives the precise
inequality. The substitutions for all five cases are in
Appendix~\ref{supp:proof:rankmanual-main}. Together with the remaining cases, they give the upper bound over
all true-null counts. The table entries are not claimed to be exact
optima for each count. The all-null vector $(U,\ldots,U)$ attains $q$, and
Proposition~\ref{prop:eight} attains $21q/20$ in dimension eight.

Thus eight is the first inflation dimension in the monotone rank-one
class. By Proposition~\ref{prop:projection}, the same upper bounds apply
when only the true-null latent covariance has rank one. They do not
settle the comparison $f_K(q)=s_K(q)$ for unrestricted Gaussian models
in dimensions four through seven.

\subsection{The sharp constant over dimensions}

Having identified the first inflation dimension in the monotone rank-one
Gaussian class, we now ask how large the inflation can become as $K$
increases. Our aim is to find the smallest bound on $\FDR/q$ that holds
uniformly over all dimensions.

Let $\omega>0$ be the solution of $\omega e^\omega=1$.
\begin{theorem}\label{thm:sharp}
For $0<q<1/2$, every model in the unrestricted-alternative two-block class
satisfies $\FDR\le 2\omega q$. This constant is sharp already in the
monotone rank-one Gaussian class:
\[
 \sup_{K\ge1}\frac{f_K^{(1)}(q)}q
 =\lim_{K\to\infty}\frac{f_K^{(1)}(q)}q=2\omega.
\]
The convergence is uniform over $0<q<1/2$. For each finite $K$, the
supremum ratio in either class is strictly less than $2\omega$.
\end{theorem}

\begin{proof}
For $m>0$, comparison with the integral of $1/t$ gives
$\sum_{j=1}^{\ell}1/(m+j)\le\log(1+\ell/m)$, strictly if $\ell>0$.
The function $t\log(1+\ell/t)$, defined to be zero at $t=0$, is
concave. Thus Theorem~\ref{thm:two-block}, with $n=a+b>0$ and $x=n/K$, gives
\[
 \frac{\FDR}{q}
 \le\frac nK\left[1+\log\left(1+\frac{2\ell}{n}\right)\right]
 =c(x):=x\left[1+\log\frac{2-x}{x}\right].
\]
Here $c(0):=0$, $c(1)=1$, and
\[
 c'(x)=\log\frac{2-x}{x}-\frac{x}{2-x},\qquad
 c''(x)=-\frac4{x(2-x)^2}<0.
\]
Setting $y=(2-x)/x$ in $c'(x)=0$ gives $y\log y=1$, so the
unique maximizer is $x_\star=2\omega/(1+\omega)$ and
$c(x_\star)=2\omega$. If $n=0$, the FDR is zero.
For fixed $K$, the finitely many integer splits with $n,\ell>0$
have a strict harmonic-sum comparison, while those with $\ell=0$
have bound one. Their maximum is therefore strictly below $2\omega$.

For sharpness, the paired construction in
Theorem~\ref{supp:thm:two-tail} takes $a$ nulls in each tail and
$m$ pairs of monotone alternatives, with $K=2a+2m$.
In either tail, its first interval of length $ah$, where $h=q/K$,
has FDP one; each subsequent interval of length $2h$ adds one
helper pair and has FDP $a/(a+2i)$. The construction keeps unmatched
helpers above the intervening BH thresholds. Hence its ratio is
\[
 L_{a,m}:=\frac{\FDR}{q}
 =\frac{2a}{K}+\frac{4a}{K}\sum_{i=1}^{m}\frac1{a+2i}.
\]
Take $K=2M\to\infty$, choose $a/M=:x_M\to x_\star$, and put
$m=M-a$. These are positive for all sufficiently large $M$, and
the Riemann sum gives
\[
 L_{a,m}=x_M\left[1+\frac2M\sum_{i=1}^{m}
 \frac1{x_M+2i/M}\right]
 \longrightarrow x_\star\left[1+\log\frac{2-x_\star}{x_\star}\right]
 =2\omega.
\]

For odd $K$, run the construction in dimension $N=K-1$ at level
$q'=qN/K$ and add one alternative identically equal to one.
For $r\le N$, the thresholds satisfy $rq/K=rq'/N$; the new
coordinate cannot be rejected. The resulting ratio is
$(N/K)L_{a,m}$ and has the same limit. The constant alternative
remains admissible in the rank-one Gaussian model. Both lower-bound
ratios are independent of $q$, although the maps are chosen for each
$q$. Together with the common upper bound, this proves the limit
uniformly over $0<q<1/2$ and hence the supremum assertion.
\end{proof}

Here $x=n/K$ is the true-null proportion. In the attaining sequence,
it approaches $x_\star=2\omega/(1+\omega)\approx0.7238$, so about
$72.4\%$ of the hypotheses are true nulls and $27.6\%$ are alternatives.
The optimization balances two effects: alternatives make larger BH
rejection ranks available, but each rejected alternative lowers the
fraction of false discoveries.

% We have shown that the global null need not be least favorable in a
% fixed dimension. The exact classification and sharp constant above
% concern the stated two-block and rank-one classes; they do not settle
% the unrestricted Gaussian comparison in dimensions four through seven.
% Our eight-dimensional separation from $s_8(q)$ is established only for
% sufficiently small $q$. Extending it to conventional testing levels
% requires a sharper global-null bound, and the general monotone
% alternatives used here do not provide a standard Gaussian mean-shift
% example. Whether a comparable separation occurs in small dimensions
% within that submodel remains a further question.

\section{Conclusion}
Under Gaussian dependence, the rejection probability of the Simes test
under the global null can understate the worst-case FDR of the BH
procedure, as recently shown by \cite{Dobriban2026} and \cite{Lei2026}. In our one-sided model with monotone alternatives, the global
null is least favorable for up to three tests. With eight tests, the
two alternatives in our construction raise the BH threshold enough to
increase false discoveries in both tails. For sufficiently small nominal
levels, the resulting FDR exceeds the largest rejection probability
attainable under the global null.

The results give benchmarks for assessing BH under dependence. At nominal
level $0.05$, the probability of any rejection under the global null is
at most $0.05147$ in any dimension when the Gaussian correlations are
nonpositive. With a fixed number of tests, our stability result gives
conditions on correlations involving true nulls under which the FDR
agrees to first order with the independence benchmark at small levels. 
When using simulations to assess BH for a particular testing design, it
is therefore important to include plausible configurations with signals.
Varying the number and strength of those signals and their dependence
with the null statistics helps assess how the BH threshold and FDR
respond to the alternatives.

\section*{Statement on AI use}
The paper is a result of a few years' thoughts and drafts (since 2020) from the authors and intensive interactions  with 
OpenAI's ChatGPT in the past year. ChatGPT helped on the technical development, editing, and mathematical checking of the paper.

\ExportMainNumbers
\appendix
\section*{Appendices}
The appendices contain the technical arguments and additional results
in the order of the main paper. The index gives the location of each proof
or supporting argument.

The general monotone stability bound is proved in
Appendix~\ref{supp:proof:thm:stability}. The stronger standard-mean bound
and its necessary-scale consequence are appendix results in
Appendices~\ref{supp:proof:thm:nondegenerate}
and~\ref{supp:proof:cor:necessary-pair}.
The density comparison in Appendix~\ref{supp:supp:sign-sectors} supports
the main-paper proof of \MainResult{thm:nonpositive-main}.
The paired construction in Appendix~\ref{supp:proof:thm:two-tail} supplies the
attaining sequence used in the main-paper proof of \MainResult{thm:sharp}.

\begingroup
\small
\setlength{\LTleft}{0pt}
\setlength{\LTright}{0pt plus 1fill}
\renewcommand{\arraystretch}{1.12}
\begin{longtable}{@{}p{0.32\linewidth}p{0.48\linewidth}r@{}}
\toprule
Main-paper result & Proof location & Page\\
\midrule
\endfirsthead
\toprule
Main-paper result & Proof location & Page\\
\midrule
\endhead
\bottomrule
\endfoot
\MainResult{enh:projection} & Proved in the main text & ---\\
\MainResult{thm:fixed} & Appendix~\ref{supp:proof:thm:fixed} & \pageref{supp:proof:thm:fixed}\\
\MainResult{enh:stability} & Appendix~\ref{supp:proof:thm:stability} & \pageref{supp:proof:thm:stability}\\
\MainResult{thm:nonpositive-main} & Main text; density comparison in Appendix~\ref{supp:supp:sign-sectors} & \pageref{supp:supp:sign-sectors}\\
\MainResult{lem:oracle} & Proved in the main text & ---\\
\MainResult{thm:one-null} & Proved in the main text & ---\\
\MainResult{thm:two-null} & Proved in the main text & ---\\
\MainResult{cor:small-K} & Proved in the main text & ---\\
\MainResult{prop:eight} & Proved in the main text & ---\\
\MainResult{cor:K8} & Proved in the main text & ---\\
\MainResult{thm:truncated} & Proved in the main text & ---\\
\MainResult{cor:balanced} & Proved in the main text & ---\\
\MainResult{rankmanual-main} & Appendix~\ref{supp:proof:rankmanual-main} & \pageref{supp:proof:rankmanual-main}\\
\MainResult{thm:sharp} & Proved in the main text & ---\\
\end{longtable}
\endgroup

% The shared proof file starts its standalone numbering at S2.
% Start the first section at A in the combined manuscript.
\let\AppendixSection\section
\renewcommand{\section}{%
  \let\section\AppendixSection
  \setcounter{section}{0}%
  \section
}
% Supplementary proofs and supporting results.
% Internal labels use supp:; main-result numbers are supplied by the entrypoint.
% This file is intended to be read with \input, not compiled by itself.
\setcounter{section}{1}
\section{Notation and global-null bounds}\label{supp:sec:model}

We recall the model and then prove the uniform global-null benchmark
used in the extremal comparison. Main-paper results retain their
main-paper numbers in references.

\subsection{Notation}
Fix $K\ge1$ and $q\in(0,1)$.  For
$p=(p_1,\ldots,p_K)\in[0,1]^K$, write
$p_{(1)}\le\cdots\le p_{(K)}$ and define
\[
  S_K(p)=\min_{1\le k\le K}\frac{Kp_{(k)}}{k}.
\]
Equivalently, with $N_p(t)=\sum_{i=1}^K\one\{p_i\le t\}$,
\[
  \{S_K(p)\le q\}
  =\left\{\exists k\in\{1,\ldots,K\}:
    N_p(kq/K)\ge k\right\}.
\]

Let $\cC_K$ be the compact set of $K\times K$ correlation matrices.
For $\Sigma\in\cC_K$, let $X\sim N(0,\Sigma)$, allowing singular
$\Sigma$, and put $P_i=\Phi(X_i)$.  Each $P_i$ is uniform; let
$\cG_K$ denote the class of all such p-value vectors.  Define the
global-null Gaussian Simes envelope
\begin{equation}\label{supp:eq:def-s}
  s_K(q)=\sup_{\Sigma\in\cC_K}
  \Pbb_\Sigma\{S_K(P)\le q\}.
\end{equation}
Independence and the arbitrary-dependence bound give
\begin{equation}\label{supp:eq:s-basic}
  q\le s_K(q)\le (H_Kq)\wedge1,
  \qquad H_K=\sum_{k=1}^K\frac1k.
\end{equation}

For the mixed-null problem, let $I_0\subseteq\{1,\ldots,K\}$ be the
true-null set and $K_0=|I_0|$. We also write $m_0=K_0$ in the
stability bounds of Section~\ref{supp:sec:nondegenerate}. Start from a centered unit-variance
Gaussian vector $X$, again allowing a singular correlation matrix.  If
$i\in I_0$, set $P_i=\Phi(X_i)$.  If $i\notin I_0$, let
$P_i=a_i(X_i)$ for an arbitrary Borel nondecreasing map
$a_i:\Rbb\to[0,1]$.  The maps and the correlation matrix may depend on
$q$.

BH rejects $R$ hypotheses, where
\[
  R=\max\{r:p_{(r)}\le rq/K\},
\]
with $R=0$ if the set is empty.  Let $F$ be the number of rejected true
nulls and set
\[
  \FDP=\frac{F}{R\vee1},\qquad \FDR=\Ebb[\FDP].
\]
Define
\begin{equation}\label{supp:eq:def-f}
  f_{K_0,K}(q)=\sup\FDR,
  \qquad
  f_K(q)=\max_{1\le K_0\le K}f_{K_0,K}(q).
\end{equation}
The supremum is over the preceding mixed Gaussian class with exactly
$K_0$ true nulls.  Under the global null, every rejection is false, so
\begin{equation}\label{supp:eq:allnull}
  f_{K,K}(q)=s_K(q),\qquad f_K(q)\ge s_K(q).
\end{equation}
The least-favorability question is whether equality always holds in the
second relation.

The standard one-sided mean-test submodel has $P_i=\barPhi(T_i)$,
$T\sim N(\mu,\Sigma)$, and boundary null means $\mu_i=0$.
Negative null means instead give superuniform margins, while two-sided
p-values are not monotone in the latent coordinate; neither is part of
the exact-uniform monotone model above.

We use a superscript $(1)$ when the centered Gaussian latent vector is
restricted to have rank one:
\[
  f^{(1)}_{K_0,K}(q)=\sup_{\operatorname{rank}(\Sigma)=1}\FDR,
  \qquad
  f^{(1)}_K(q)=\max_{1\le K_0\le K}f^{(1)}_{K_0,K}(q).
\]
Every standardized coordinate of such a vector is either $Z$ or $-Z$
for one standard normal $Z$. This two-tail representation supports the
analytical classification and sharp bounds in Section~\ref{supp:sec:rankone}.

The projection reduction and its proof are given in
\MainResult{enh:projection} in the main paper.

\subsection{Proof of \MainResult{thm:fixed}}\label{supp:proof:thm:fixed}\label{supp:sec:fixed}

The fixed-matrix small-level Simes upper bound is due to \citet[Theorem~8(III)]{ChiRamdasWang2024}.
Here we give the uniform remainder needed for the global-null
comparison in the main paper. Grouping nearly perfectly positively
correlated coordinates gives an error bound depending only on
$K$ and $q$, even when the matrix varies with $q$.

\MainResult{thm:fixed} follows from the following uniform estimate:
its error is $o(q)$ for fixed $K$, and independence gives the matching
lower bound. The case $K=1$ is immediate.

\begin{theorem}\label{supp:ms:thm:global}
For every integer $K\ge2$, let
$c_K=\sqrt{2/(2-(K-1)^{-2})}$. If $0<q\le1/2$ and
$t=\barPhi^{-1}(q)$, then
\begin{equation}\label{supp:ms:eq:global-tail}
 s_K(q)\le q+\frac{K(K-1)}2\barPhi(c_Kt).
\end{equation}
For $t\ge1$, in particular,
\begin{equation}\label{supp:ms:eq:global-rate}
 \frac{s_K(q)}q
 \le1+K(K-1)\exp\left\{-\frac{t^2}{4(K-1)^2-2}\right\}.
\end{equation}
Consequently $s_K(q)<1.001q$ whenever
$t\ge2K\sqrt{\log(1000K^2)}$.
\end{theorem}

\begin{proof}
We may use upper-tail p-values $\barPhi(X_i)$ in place of
$\Phi(X_i)$: replacing the entire centered vector by its negative
does not change its correlation matrix or Gaussian law.
Represent any correlation matrix as the Gram matrix of unit vectors
$v_1,\ldots,v_K$. Join $i,j$ when
$\rho_{ij}>1-(K-1)^{-2}$. Each edge has Euclidean length
less than $\sqrt2/(K-1)$. Two vertices in the same connected component
have a simple path with at most $K-1$ edges, so their distance is
less than $\sqrt2$ and their inner product is positive. Different
components have correlations at most $1-(K-1)^{-2}$.

If a component has $k_g$ coordinates, its Gaussian correlations
are nonnegative. The positive-dependence Simes inequality
\citep{BY2001}, applied at level $qk_g/K$, bounds the probability of its event
$\{\exists r:P^{(g)}_{(r)}\le rq/K\}$ by $qk_g/K$.
The same bound holds for singular matrices by adding independent
Gaussian noise and taking a limit; the boundary consists of finitely
many marginal grid levels, each with probability zero.

A full Simes crossing either lies within one component, or contains
two coordinates from different components whose p-values are at most
$q$. For such a pair, the sum of their statistics gives
\[
 \Pbb(X_i\ge t,X_j\ge t)
 \le\barPhi\left(\sqrt{2/(1+\rho_{ij})}\,t\right)
 \le\barPhi(c_Kt).
\]
When $\rho_{ij}=-1$, the pair probability is zero. Summing component
bounds and at most $K(K-1)/2$ pair bounds proves
\eqref{supp:ms:eq:global-tail} uniformly over every correlation matrix.

The Mills bounds
\begin{equation}\label{supp:ms:eq:basic-mills}
 \frac{x}{1+x^2}\phi(x)<\barPhi(x)<\frac{\phi(x)}x,
 \qquad x>0,
\end{equation}
give $\barPhi(c_Kt)/\barPhi(t)\le
2\exp\{-(c_K^2-1)t^2/2\}$ when $t\ge1$.
Since $(c_K^2-1)/2=1/[4(K-1)^2-2]$, this proves
\eqref{supp:ms:eq:global-rate}. At the stated threshold, the excess over one is less
than $K(K-1)/(1000K^2)<1/1000$.
\end{proof}

This estimate provides the uniform remainder used in the global-null
comparison. No sharpness in $K$ or $q$ is claimed, and the estimate
need not improve on the Hommel--BY bound at conventional testing levels.

\subsection{Further global-null reductions}\label{supp:sec:globalstructure}

Repeating p-values preserves Gaussian dependence and changes the
weights with which distinct coordinates enter the Simes statistic.
The following reductions explain how a fixed-dimensional example
can be embedded in larger dimensions. They also give a weighted
description of the supremum over dimensions.

\subsubsection{Clone and padding reductions}

\begin{lemma}\label{supp:lem:clone}
If $m$ divides $K$, then
\[
  s_K(q)\ge s_m(q).
\]
More precisely, repeating every coordinate of an $m$-dimensional
Gaussian-copula vector $K/m$ times leaves its Simes statistic
unchanged.
\end{lemma}

\begin{proof}
Write $K=rm$.  If $U_{(1)}\le\cdots\le U_{(m)}$, the cloned order
statistics satisfy
\[
  P_{(r(j-1)+1)}=\cdots=P_{(rj)}=U_{(j)}.
\]
For a fixed $j$, the smallest value of $KP_{(k)}/k$ over this block is
attained at $k=rj$.  Therefore
\[
  S_K(P)=\min_{1\le j\le m}\frac{K U_{(j)}}{rj}
  =S_m(U).
\]
The cloned vector is a valid, possibly singular, Gaussian copula.
\end{proof}

\begin{lemma}\label{supp:lem:padding}
For $1\le m\le K$, with $r=\lfloor K/m\rfloor$,
\begin{equation}\label{supp:eq:padding}
  s_K(q)\ge s_m\left(\frac{rm}{K}q\right)
  \ge s_m\left(\left(1-\frac{m-1}{K}\right)q\right).
\end{equation}
\end{lemma}

\begin{proof}
Repeat each coordinate of an $m$-dimensional vector $r$ times and fill
the remaining $K-rm$ positions with arbitrary Gaussian-copula
coordinates.  If $U_{(j)}\le jr q/K$, then at least $rj$ of the
$K$ p-values are at most $(rj)q/K$, so the $K$-dimensional Simes event
occurs.  Finally,
$rm=K-(K\bmod m)\ge K-m+1$ and $s_m$ is nondecreasing.
\end{proof}

\begin{theorem}\label{supp:thm:limit}
For every fixed $q\in(0,1)$,
\[
  \lim_{K\to\infty}s_K(q)=\sup_{m\ge1}s_m(q)=:s_\infty(q),
\]
where convergence holds along every integer sequence $K\to\infty$.
\end{theorem}

\begin{proof}
Let $L(q)=\sup_m s_m(q)$.  By definition,
$\limsup_Ks_K(q)\le L(q)$.  Fix $\varepsilon>0$ and choose a finite
$m$ and a Gaussian copula $U$ such that
\[
  \Pbb\{S_m(U)\le q\}>L(q)-\varepsilon.
\]
The padding construction gives, for
$r_K=\lfloor K/m\rfloor$,
\[
  s_K(q)\ge
  \Pbb\left\{S_m(U)\le\frac{r_Km}{K}q\right\}.
\]
The multiplier tends to one.  Moreover, $S_m(U)$ has no atom at $q$:
\[
  \{S_m(U)=q\}
  \subseteq
  \bigcup_{i=1}^m\bigcup_{j=1}^m
  \{U_i=jq/m\},
\]
and every marginal is continuous.  The last probability therefore
converges to $\Pbb\{S_m(U)\le q\}$.  Hence
$\liminf_Ks_K(q)\ge L(q)-\varepsilon$. Letting
$\varepsilon\downarrow0$ completes the proof.
\end{proof}

\subsubsection{Weighted clone blocks}

For $m\ge1$ and a rational vector
$w=(w_1,\ldots,w_m)$ in the open probability simplex, write
$W_A=\sum_{i\in A}w_i$ and define
\begin{equation}\label{supp:eq:sigma-m}
\begin{split}
  \sigma_m(q)=\sup_w\sup_{U\in\cG_m}
  \Pbb\bigl\{\exists\,\varnothing\ne A\subseteq\{1,\ldots,m\}:\ 
  U_j\le qW_A\ \text{for every }j\in A\bigr\}.
\end{split}
\end{equation}

\begin{theorem}\label{supp:thm:weighted}
For every $q\in(0,1)$,
\[
  s_\infty(q)=\sup_{m\ge1}\sigma_m(q).
\]
Moreover, $\sigma_m(q)$ is nondecreasing in $m$.
\end{theorem}

\begin{proof}
Taking $m=K$ and equal weights $w_i=1/K$ makes the event in
\eqref{supp:eq:sigma-m} exactly the $K$-dimensional Simes event.  Thus
$s_\infty(q)\le\sup_m\sigma_m(q)$.

Conversely, write $w_i=r_i/R$ with positive integers $r_i$ and
$R=\sum_i r_i$.  Replace $U_i$ by $r_i$ clones.  If the weighted event
holds for a set $A$, then at least $\sum_{i\in A}r_i=RW_A$ p-values are
at most $qW_A$, which is precisely an $R$-dimensional Simes crossing.
Therefore $\sigma_m(q)\le s_\infty(q)$.

For monotonicity, fix an $m$-coordinate configuration and append an
independent Gaussian coordinate.  For rational
$\delta\downarrow0$, give the new coordinate weight $\delta$ and
rescale the old weights to $(1-\delta)w_i$.  Ignoring the new
coordinate recovers the old weighted event with $q$ replaced by
$q(1-\delta)$.  This event is a finite union of rectangles whose
boundary has probability zero, so its probability converges to that
of the original event.  Taking suprema proves
$\sigma_{m+1}(q)\ge\sigma_m(q)$.
\end{proof}

\section{Stability away from perfect correlation}
\label{supp:nd:sec:proof}\label{supp:sec:nondegenerate}

We first establish a sharper stability bound for the standard Gaussian
mean-shift subclass and derive the necessary correlation scale. We then
prove the main-paper result \MainResult{enh:stability} for general monotone
margins.

The common idea is to replace a true null's random effect on the
rejection count by its value at the tail threshold. The resulting
reference contribution is exactly $q/K$. We establish the forced-zero
identity, bound the probability that freezing changes the count, and
integrate that error.

\subsection{A sharper bound for standard Gaussian mean tests}\label{supp:proof:thm:nondegenerate}

Let $K\ge2$, $T\sim N(\mu,\Sigma)$, and $P_i=\barPhi(T_i)$, with
$\mu_i=0$ for $i\in I_0$ and $\Sigma_{ii}=1$. The full covariance
matrix and the false-null subvector may be singular.
Let $I_0$ be fixed, $m_0=|I_0|$, and suppose
$|\rho_{ij}|<1$ for $i\in I_0$, $j\ne i$.  Put
\[
 t=\barPhi^{-1}(q),\quad M(q)=\phi(t)-tq,
 \quad w_K(q)=\barPhi^{-1}(2q/K)-t.
\]
Define
\begin{align*}
 D_K(\Sigma,I_0)
 &=\frac1{\sqrt{2\pi}}
 \sum_{i\in I_0}\sum_{j\ne i}
 \frac{\min\{(K-1)|\rho_{ij}|,
       |\rho_{ij}|+\log(K/2)\}}
      {\sqrt{1-\rho_{ij}^2}},\\
 C_{K,m_0}
 &=\frac{m_0(K-1)\{1+\log(K/2)\}}{\sqrt{2\pi}}.
\end{align*}

\begin{theorem}
\label{supp:nd:thm:main}
For $0<q<1/2$,
\begin{align*}
 \left|\FDR-\frac{m_0q}{K}\right|
 &\le \frac1{\sqrt{2\pi}}
 \sum_{i\in I_0}\sum_{j\ne i}
 \frac{\min\{(K-1)|\rho_{ij}|M(q),
       |\rho_{ij}|M(q)+qw_K(q)\}}
      {\sqrt{1-\rho_{ij}^2}}\\
 &\le\frac qtD_K(\Sigma,I_0).
\end{align*}
The bounds are uniform over all alternative means, including means
depending on $q$.
\end{theorem}

\begin{lemma}
\label{supp:nd:lem:loo}
For a deterministic p-value vector, let $R_i^0$ be the BH rejection count
after replacing $p_i$ by zero.  Then $R_i^0\ge1$ and
\begin{equation}\label{supp:nd:eq:1}
 \frac{\one\{i\text{ is rejected}\}}{R\vee1}
 =\frac{\one\{p_i\le qR_i^0/K\}}{R_i^0}.
\end{equation}
The count $R_i^0$ is determined by the comparisons
$p_j\le rq/K$, $j\ne i$, $r=2,\ldots,K$.
\end{lemma}

\begin{proof}
Lowering one p-value cannot decrease the step-up rejection count.  If $i$
is rejected in the original vector, then replacing $p_i$ by zero creates no
new crossing above the original rejection rank, because $p_i$ was already
counted at every larger threshold.  Hence $R_i^0=R$.  Conversely, if
$p_i\le qR_i^0/K$, restoring $p_i$ preserves the crossing at rank $R_i^0$,
so monotonicity gives $R=R_i^0$ and $i$ is rejected.  A crossing at rank
$r\ge2$ in the forced-zero vector occurs precisely when at least $r-1$ of
the other p-values do not exceed $rq/K$.
\end{proof}

Fix $i\in I_0$.  Gaussian regression gives a vector $Y^{(i)}$, independent
of $T_i$, such that
\begin{equation}\label{supp:nd:eq:2}
 T_j=\mu_j+\rho_{ij}T_i+Y_j^{(i)},\qquad
 \operatorname{Var}(Y_j^{(i)})=1-\rho_{ij}^2.
\end{equation}
The coordinates of $Y^{(i)}$ need not be mutually independent.
Let $\mathcal R_i(x,Y^{(i)})$ be the forced-zero rejection count when the
other statistics equal $\mu_j+\rho_{ij}x+Y_j^{(i)}$.

\begin{lemma}
\label{supp:nd:lem:freeze}
For every $x\ge t$,
\begin{equation}\label{supp:nd:eq:3}
\begin{aligned}
 &\Pbb\{\mathcal R_i(x,Y^{(i)})\ne
              \mathcal R_i(t,Y^{(i)})\}\\
 &\quad\le\frac1{\sqrt{2\pi}}
 \sum_{j\ne i}
 \frac{\min\{(K-1)|\rho_{ij}|(x-t),
       |\rho_{ij}|(x-t)+w_K(q)\}}
      {\sqrt{1-\rho_{ij}^2}}.
\end{aligned}
\end{equation}
\end{lemma}

\begin{proof}
For $r=2,\ldots,K$, put $c_r=\barPhi^{-1}(rq/K)$.  If the forced-zero
counts differ, at least one comparison
\[
 \mu_j+\rho_{ij}x+Y_j^{(i)}\ge c_r
\]
must change between $x$ and $t$.  For fixed $(j,r)$, the residual must lie in
an interval of length $|\rho_{ij}|(x-t)$.  Its one-dimensional density is
bounded by $\{2\pi(1-\rho_{ij}^2)\}^{-1/2}$.  Summing over the $K-1$ ranks
gives the first argument of the minimum.

For a fixed $j$, the union of all threshold intervals lies between the two
translated endpoint thresholds.  Its total length is at most
$|\rho_{ij}|(x-t)+c_2-c_K$.  Since $c_K=t$ and
$c_2=t+w_K(q)$, this gives the second argument.  Taking the better of the
two bounds and summing over $j$ proves the lemma.  The interval locations
depend on $\mu_j$, but their lengths and density bounds do not.  The bound
is therefore uniform over the means.
\end{proof}

\begin{lemma}
\label{supp:nd:lem:tails}
For $0<q<1/2$,
\begin{equation}\label{supp:nd:eq:4}
 \int_t^\infty(x-t)\phi(x)\,dx=M(q),
 \qquad 0<M(q)<\frac qt,
\end{equation}
and
\begin{equation}\label{supp:nd:eq:5}
 0\le w_K(q)\le\frac{\log(K/2)}t.
\end{equation}
\end{lemma}

\begin{proof}
Integration by parts gives the first identity.  The strict Mills lower
bound $q>t\phi(t)/(1+t^2)$ yields $M(q)<q/t$.  To bound the width, let
$h(x)=\phi(x)/\barPhi(x)$ be the normal hazard rate.  Since $h(x)\ge x$ for
$x>0$,
\[
 \log(K/2)=\log\frac{\barPhi(t)}{\barPhi(t+w_K)}
 =\int_t^{t+w_K}h(x)\,dx\ge tw_K.
\]
The case $K=2$ has $w_K=0$.
\end{proof}

\begin{proof}[Proof of Theorem~\ref{supp:nd:thm:main}]
By Lemma~\ref{supp:nd:lem:loo}, the contribution of true null $i$ is
\[
 \Ebb\left[
 \frac{\one\{T_i\ge\barPhi^{-1}(qR_i^0/K)\}}{R_i^0}
 \right].
\]
For comparison, replace
$R_i^0=\mathcal R_i(T_i,Y^{(i)})$ with the frozen count
$\mathcal R_i(t,Y^{(i)})$.  Conditional on $Y^{(i)}$, this count is an
integer $r\in\{1,\ldots,K\}$ independent of $T_i$, so the conditional
expectation of the reference contribution is
\begin{equation}\label{supp:nd:eq:6}
 \frac1r\Pbb\{T_i\ge\barPhi^{-1}(qr/K)\}=\frac qK.
\end{equation}
Both the true and frozen integrands vanish when $T_i<t$.  Their difference
is at most the indicator that the two counts differ.  Lemma
\ref{supp:nd:lem:freeze}, followed by integration against $\phi(x)$ on
$[t,\infty)$, gives
\begin{align*}
 \left|\text{contribution of }i-\frac qK\right|
 &\le\frac1{\sqrt{2\pi}}
 \sum_{j\ne i}\frac1{\sqrt{1-\rho_{ij}^2}}\\
 &\quad\times
 \min\{(K-1)|\rho_{ij}|M(q),
       |\rho_{ij}|M(q)+qw_K(q)\}.
\end{align*}
Summing over $i\in I_0$ proves the finite-level bound.  Lemma
\ref{supp:nd:lem:tails} gives the displayed rate for this standard mean-shift model.
\end{proof}

\subsection{The necessary correlation scale}\label{supp:proof:cor:necessary-pair}

\begin{corollary}
\label{supp:nd:cor:necessary}
Let
$d(\Sigma,I_0)=\min_{i\in I_0,j\ne i}(1-\rho_{ij}^2)$.  If
\[
 \frac{\FDR}{q}\ge\frac{m_0}{K}+\delta,
 \qquad\delta>0,
\]
then
\[
 d(\Sigma,I_0)\le
 \min\left\{1,\frac{C_{K,m_0}^2}{\delta^2t^2}\right\}.
\]
Along any fixed-dimensional family with a fixed positive excess, some pair
involving a true null therefore has
$1-\rho_{ij}^2=O\{1/\log(1/q)\}$, and
$\lambda_{\min}(\Sigma)=O\{1/\log(1/q)\}$.
\end{corollary}

\begin{proof}[Proof of Corollary~\ref{supp:nd:cor:necessary}]
Let
\[
 d=\min_{i\in I_0,\,j\ne i}(1-\rho_{ij}^2).
\]
If $d=0$, the conclusion
is immediate.  Otherwise,
\[
 D_K(\Sigma,I_0)\le
 \frac{m_0(K-1)\{1+\log(K/2)\}}{\sqrt{2\pi d}}
 =\frac{C_{K,m_0}}{\sqrt d}.
\]
The assumed excess and the theorem imply
$\delta\le C_{K,m_0}/(t\sqrt d)$, which proves the finite bound.  Since
$t^2\sim2\log(1/q)$, the asserted asymptotic scale follows.  For a pair
achieving $d$, a signed two-coordinate unit vector gives
\[
 \lambda_{\min}(\Sigma)\le1-|\rho_{ij}|
 \le1-\rho_{ij}^2=d.
\]
There are finitely many pairs, so one pair is selected along a subsequence.
\end{proof}

\begin{corollary}
For fixed $K,I_0$, suppose that
\[
 d(\Sigma(q),I_0)\log(1/q)\longrightarrow\infty.
\]
Then
$\FDR/q\to m_0/K$, uniformly over the alternative means.  If $m_0<K$,
eventually $\FDR<q$.
\end{corollary}

\begin{proof}
The condition gives $t\sqrt d\to\infty$, so the rate bound above applies.
\end{proof}
\subsection{Proof of \MainResult{enh:stability}}\label{supp:ed:stability-proof}\label{supp:proof:thm:stability}

Use $T=-X$ so true-null p-values are $\barPhi(T_i)$ and alternative
maps $b_j(T_j)=a_j(-T_j)$ are nonincreasing. For each BH threshold
$rq/K$, every event $\{b_j(T_j)\le rq/K\}$ has an upper half-line
preimage, possibly empty or all of $\mathbb R$. Denote its endpoint
by $c_{jr}\in[-\infty,\infty]$. Endpoint conventions do not matter
for the Gaussian comparisons used below. True-null coordinates have
the same representation with their corresponding normal quantiles.

Fix $i\in I_0$ and regress
\[
 T_j=\rho_{ij}T_i+Y_j,\quad j\ne i.
\]
The residual vector is independent of $T_i$ and its $j$th marginal
variance is $1-\rho_{ij}^2>0$. Its joint law need not be nonsingular.
Let $R_i^0(x,Y)$ be the rejection count with coordinate $i$ forced to
zero and all other latent values evaluated at $\rho_{ij}x+Y_j$.
This count is at least one and depends only on comparisons at ranks
$r=2,\ldots,K$.

Write $t=\barPhi^{-1}(q)>0$. For $x\ge t$, a change between
$R_i^0(x,Y)$ and $R_i^0(t,Y)$ requires at least one comparison to
cross its endpoint. For a finite $c_{jr}$, the necessary residual
interval has length $|\rho_{ij}|(x-t)$. Its probability is at most
this length divided by $\sqrt{2\pi(1-\rho_{ij}^2)}$.
An infinite endpoint cannot change. The union bound therefore gives,
without assuming residual independence,
\begin{equation}\label{supp:enh:supp-freeze}
 \Pbb\{R_i^0(x,Y)\ne R_i^0(t,Y)\}
 \le\frac{K-1}{\sqrt{2\pi}}
 \sum_{j\ne i}\frac{|\rho_{ij}|(x-t)}{\sqrt{1-\rho_{ij}^2}}.
\end{equation}

Lemma~\ref{supp:nd:lem:loo} gives the actual
FDR contribution of null $i$ as
\[
 \Ebb\left[
 \frac{\one\{\barPhi(T_i)\le qR_i^0(T_i,Y)/K\}}
      {R_i^0(T_i,Y)}\right].
\]
Replace that count by $R_i^0(t,Y)$. Conditional on $Y$, the frozen
count $r$ is independent of $T_i$, and the frozen expectation is
$(qr/K)/r=q/K$. Both integrands vanish for $T_i<t$, and their
absolute difference is at most the indicator that the counts differ.
Integrating \eqref{supp:enh:supp-freeze} over $x\ge t$ gives
\[
 \left|\FDR-\frac{m_0q}{K}\right|
 \le\frac{(K-1)M(q)}{\sqrt{2\pi}}
 \sum_{i\in I_0}\sum_{j\ne i}
 \frac{|\rho_{ij}|}{\sqrt{1-\rho_{ij}^2}},
 \quad M(q)=\int_t^\infty(x-t)\phi(x)\dd x=\phi(t)-tq.
\]
The Mills lower bound $q>t\phi(t)/(1+t^2)$ implies $M(q)<q/t$.
Bounding the $m_0(K-1)$ ordered summands by $1/\sqrt d$ proves
\[
 \left|\frac{\FDR}{q}-\frac{m_0}{K}\right|
 \le\frac{m_0(K-1)^2}{\sqrt{2\pi}\,t\sqrt d},
 \quad d=\min_{i\in I_0,j\ne i}(1-\rho_{ij}^2).
\]
All estimates are independent of the maps' cutpoint positions, so the
maps may depend on $q$ and may have jumps, atoms or flat regions.
For $d\ge\eta(q)$, the last bound is at most
$m_0(K-1)^2/\{\sqrt{2\pi}\,t\sqrt{\eta(q)}\}$. This tends to zero
when $\eta(q)\log(1/q)\to\infty$, since $t^2\sim2\log(1/q)$,
and proves the uniform convergence in \MainResult{enh:stability}.
The necessary scale follows by rearrangement for a fixed positive
excess above $m_0/K$. When $d=0$ it is automatic. This argument does
not assert the common-width improvement available for standard shifts.

\section{Gaussian density comparison}\label{supp:supp:sign-sectors}\label{supp:proof:thm:nonpositive-main}

The proof of \MainResult{thm:nonpositive-main} is given in the main text.
We establish here the density comparison used to bound its short-interval
probabilities by volumes. Every principal submatrix inherits the
assumptions of the following lemma.

\begin{lemma}\label{supp:lem:density-dom}
Let $\Sigma$ be a positive-definite $d\times d$ correlation matrix
with nonpositive off-diagonal entries, and let $\varphi_\Sigma$ be its
Gaussian density. If $x_i\le-1$ for every $i$, then
\begin{equation}\label{supp:eq:density-dom}
 \varphi_\Sigma(x)\le\prod_{i=1}^d\varphi(x_i).
\end{equation}
\end{lemma}

\begin{proof}
Let $R=\Sigma-I$, $\Sigma_t=I+tR$, and
$B_t=\Sigma_t^{-1}$.  The positive-definite matrix $\Sigma_t$ has
nonpositive off-diagonal entries, so it is a Stieltjes matrix and
$B_t$ is entrywise nonnegative.  Also $B_t\Sigma_t=I$ implies
$(B_t)_{ii}\ge1$.  Put $y=-x\ge\mathbf1$ and
\[
  L(t)=\log\frac{\varphi_{\Sigma_t}(x)}{\prod_i\varphi(x_i)}.
\]
Differentiation gives
\[
  L'(t)=\frac12\sum_{i,j}R_{ij}
  \left\{(B_ty)_i(B_ty)_j-(B_t)_{ij}\right\}.
\]
The diagonal terms vanish.  For $i\ne j$,
\[
  (B_ty)_i(B_ty)_j
  \ge (B_t)_{ij}(B_t)_{jj}\ge(B_t)_{ij},
\]
whereas $R_{ij}\le0$.  Thus $L'(t)\le0$ and
$L(1)\le L(0)=0$, proving \eqref{supp:eq:density-dom}.

\end{proof}

\section{Mixed-null reductions}\label{supp:sec:mixed}

The oracle bound and the one- and two-null reductions are proved in
\MainResult{lem:oracle}, \MainResult{thm:one-null},
\MainResult{thm:two-null}, and \MainResult{cor:small-K} in the main paper.
This section gives a related integral bound for arbitrary true-null counts.
Let $n$ be the number of true nulls among the $K$ hypotheses, and let
$Q$ be their p-value vector, with ordered coordinates
$Q_{(1)}\le\cdots\le Q_{(n)}$.

An expectation can also be controlled by asking how likely the FDP
is to exceed each value $x$. A large FDP forces a Simes crossing
among the true nulls at a suitable level. Integrating this implication
produces a bound in terms of $s_n$ alone.

The following supporting bound implements the null-to-full FDR-linking
principle \citep{Su2018} with finite counts. It does not introduce a new
general linking method.

For $x\in(0,1]$ and $j=1,\ldots,n$, define
\[
  r_{j,n,K}(x)=
  \min\left\{j+K-n,\left\lfloor\frac{j}{x}\right\rfloor\right\}
\]
and
\[
  u_{n,K,q}(x)=
  \left[
  \frac{nq}{K}\max_{1\le j\le n}\frac{r_{j,n,K}(x)}{j}
  \right]\wedge1.
\]
Adopt the convention $s_n(1)=1$.

\begin{proposition}\label{supp:prop:layercake}
For $1\le n\le K$,
\[
  f_{n,K}(q)\le\int_0^1s_n\{u_{n,K,q}(x)\}\dd x.
\]
\end{proposition}

\begin{proof}
Let $Y=F/(R\vee1)$.  If $Y\ge x>0$, then for
$j=F\ge1$ and $r=R$,
\[
  r-j\le K-n,\qquad r\le j/x,
\]
and therefore $r\le r_{j,n,K}(x)$.  The $j$th ordered true-null
p-value obeys $Q_{(j)}\le qr/K$.  It follows that
\[
  S_n(Q)\le\frac{nQ_{(j)}}j
  \le\frac{nq}{K}\frac{r_{j,n,K}(x)}j
  \le A(x):=\frac{nq}{K}\max_{1\le k\le n}
                 \frac{r_{k,n,K}(x)}{k}.
\]
Since also $S_n(Q)\le1$, this yields
$S_n(Q)\le A(x)\wedge1=u_{n,K,q}(x)$.
Thus $\Pbb(Y\ge x)\le s_n\{u_{n,K,q}(x)\}$. The tail-probability identity
$\Ebb Y=\int_0^1\Pbb(Y\ge x)\dd x$ gives the result.
\end{proof}

\section{Rank-one constructions and bounds}
\label{supp:sec:twotail}\label{supp:sec:rankone}

We first identify the limitation of a single tied null block. We then
use two disjoint tails, where helper pairs extend the rejection region
while preserving a sufficiently large fraction of false discoveries.
The dimension-eight example is proved in \MainResult{prop:eight} in the
main paper. The paired construction supplies the sequence used to
attain the sharp constant over dimensions. We then prove the
small-dimensional classification in Section~\ref{supp:proof:rankmanual-main}.

A single comonotone null block does not create inflation: when $m$
true-null p-values all equal a uniform $U$, the following
rank-by-rank bound is sharp.

\subsection{A single tied null block}\label{supp:proof:prop:staircase}

\begin{proposition}\label{supp:prop:staircase}
If the $m$ true-null p-values are all equal to one uniform variable,
then, for arbitrary false-null p-values,
\begin{equation}\label{supp:eq:staircase}
  \FDR\le
  \frac{mq}{K}\{1+H_K-H_m\}\le q
\end{equation}
and equality in the first inequality is attainable in the monotone
Gaussian class.
\end{proposition}

\begin{proof}
Put $h=q/K$.  If the common null p-value
$U\in((r-1)h,rh]$ is rejected, then all $m$ true nulls are rejected
and $R\ge\max(m,r)$.  Hence
\[
  \FDR\le\sum_{r=1}^K
  \frac{m}{\max(m,r)}h
  =\frac{mq}{K}\{1+H_K-H_m\}.
\]
The last expression is at most $q$, since
$H_K-H_m\le\log(K/m)$ and
$x\{1+\log(1/x)\}\le1$ for $0<x\le1$.

For sharpness, write $U=\Phi(Z)$ and, for
$j=1,\ldots,K-m$, use the false-null helper
\[
  A_j(U)=
  \begin{cases}
  0,&U>(m+j-1)h,\\
  1,&U\le(m+j-1)h.
  \end{cases}
\]
It is nondecreasing in the Gaussian coordinate $-Z$.  On the $r$th
interval, BH rejects exactly $m$ hypotheses if $r\le m$ and exactly
$r$ hypotheses if $r>m$.  Thus the displayed upper bound is attained
interval by interval.
\end{proof}

With two countermonotone null blocks, helpers activated from opposite
directions can raise the rejection rank while keeping a high proportion
of true nulls among the rejections.

\subsection{Paired monotone alternatives}\label{supp:proof:thm:two-tail}

\begin{theorem}\label{supp:thm:two-tail}
Let $a,m$ be positive integers,
\[
  K=2a+2m,\qquad K_0=2a,
\]
and let $0<q<1/2$.  There is an admissible rank-one Gaussian model
whose FDR is exactly
\begin{equation}\label{supp:eq:two-tail-finite}
  \frac{\FDR}{q}
  =\frac{2a}{K}+\frac{4a}{K}
  \sum_{i=1}^m\frac1{a+2i}.
\end{equation}
\end{theorem}

\begin{proof}
Let $Z\sim N(0,1)$, $U=\Phi(Z)$, $h=q/K$, and take the true-null
p-values to be $a$ copies of $U$ and $a$ copies of $1-U=\Phi(-Z)$.
Fix $0<\varepsilon<h$.  For $i=1,\ldots,m$, define
\[
  \gamma_i=(a+2i-2)h,
  \qquad
  \lambda_i=(a+2i)h-\varepsilon,
\]
and the two false-null p-values
\[
 A_i^+(U)=
 \begin{cases}
  1,&U\le\gamma_i,\\
  \lambda_i,&U>\gamma_i,
 \end{cases}
 \qquad
 A_i^-(U)=
 \begin{cases}
  \lambda_i,&U<1-\gamma_i,\\
  1,&U\ge1-\gamma_i.
 \end{cases}
\]
Here $A_i^+$ is a nondecreasing function of the Gaussian coordinate
$-Z$, while $A_i^-$ is a nondecreasing function of $Z$.  Thus all
coordinates belong to the model of Section~\ref{supp:sec:model}.

Consider the left tail.  On $0<U\le ah$, BH rejects exactly the $a$
true nulls.  Indeed, the $s$th active unmatched helper has value
$\lambda_s=(a+2s)h-\varepsilon>(a+s)h$ for every $s\ge1$, so no later
rank crosses.
For
\[
  (a+2i-2)h<U\le(a+2i)h,
  \qquad i=1,\ldots,m,
\]
the $a$ left-tail nulls and the first $i$ helper pairs pass at rank
$a+2i$.  Thus BH makes $a+2i$ rejections and the FDP is
$a/(a+2i)$.  No skipped-rank crossing is possible: the $s$th
unmatched helper after this layer has value
\[
  \lambda_{i+s}=(a+2i+2s)h-\varepsilon,
\]
which exceeds the corresponding possible intervening critical value
$(a+2i+s)h$ by $sh-\varepsilon>0$.  After the final layer, any
helper-only rejections have FDP zero.  To see that these intervals
exhaust all possible left-tail false discoveries, observe that the two
null blocks cannot both be rejected when $q<1/2$.  Thus rejection of
the $U$ block entails at most $K-a=a+2m$ total rejections and hence
$U\le(K-a)h$, exactly the endpoint of the last displayed layer.
Here $U<(K-a)q/K<q<1/2$ and every $\gamma_i<q$, so all required
$A_i^-$ helpers are active, while $1-U>q$ and the opposite null block
and inactive helpers cannot participate.

The right tail is symmetric and disjoint.  Each tail therefore
contributes a base interval of length $ah$ with FDP one and, for each
$i$, an interval of length $2h$ with FDP $a/(a+2i)$.  Consequently
\[
  \FDR
  =2\left\{ah+\sum_{i=1}^m2h\frac{a}{a+2i}\right\},
\]
which is \eqref{supp:eq:two-tail-finite}.
\end{proof}

\subsection{Proof of \MainResult{rankmanual-main}}
\label{supp:proof:rankmanual-main}
\label{supp:sec:rankmanual}

We prove \MainResult{rankmanual-main} using bounds for each tail and
one refinement that links the two tails. The proof is analytical and
uses neither a finite-state computation nor a classification for every
fixed true-null count.

Throughout this subsection, $0<q<1/2$, $h=q/K$, and the model is the
monotone rank-one Gaussian class. There are $a$ true-null p-values
equal to $U$, $b$ equal to $1-U$, and $\ell=K-a-b$ false-null p-values,
each a nondecreasing or nonincreasing Borel function of $U$.

We first establish tail bounds.  In one tail $0<v\le q$, write
$v=U$ or $v=1-U$
for its distance coordinate, $x=v/h$, and $m\ge1$ for the number of
true nulls in that block.  The opposite block exceeds $q$ throughout
this tail and can be replaced by p-values equal to one.  The total
dimension remains $K$, with the same critical values $rh$.
Tail integrals are over $0<x\le K$; the true-null contribution
vanishes for $x>m+\ell$.

The variable $x$ expresses the common null p-value in units of the
BH grid width $h$. The function $r(x)$ below is a rejection count
after the tied null block is set to zero; it is not a p-value
threshold. Restoring the block succeeds precisely when its scaled
value $x$ is no larger than that count. When $r(x)$ is nonincreasing,
the endpoint $T$ of this success interval controls the entire
integral, even if the helper margins have jumps.

\begin{lemma}
\label{supp:rankmanual-rank}
At a fixed $x$, replace the $m$ true-null p-values in the block by
zeros, and let $r(x)$ be the resulting BH rejection count.  Then
$m\le r(x)\le m+\ell$, and the contribution of this block to the
actual FDP is exactly
\begin{equation}\label{supp:rankmanual-effective}
 \psi(x)=\frac{m}{r(x)}\one\{x\le r(x)\}.
\end{equation}
If all false-null p-values are nondecreasing in $x$, the tail FDR
is at most $mh$.
\end{lemma}

\begin{proof}
If $x\le r(x)$, restoring the true-null values to $hx$ preserves
the crossing at rank $r(x)$.  Increasing p-values cannot increase
the BH rejection count, so the actual rejection count is $r(x)$.
If $x>r(x)$, rejection of a true null would require an actual count
$R\ge x>r(x)$, which is impossible.  The tied null block is either
entirely rejected or entirely retained, proving
\eqref{supp:rankmanual-effective}.

In the monotone case, $r(x)$ is nonincreasing.  Its success set
$\{x:x\le r(x)\}$ is an interval from zero, possibly excluding its
right endpoint.  Let $T$ be that endpoint.  Then
$m\le T\le m+\ell$, and $r(x)\ge T$ for every $x<T$:
indeed, for $x<y<T$, monotonicity gives $r(x)\ge r(y)\ge y$,
and one may let $y\uparrow T$.  Thus
\[
 \int\psi(x)\,dx=\int_0^T\frac{m}{r(x)}\,dx\le m.
\]
Multiplication by $h$ gives the bound on the tail FDR. No continuity is needed;
the open or closed status of the single endpoint has no effect on
the integral.
\end{proof}

Recall the truncated staircase quantity $G_\ell$ and define $P_\ell$ by
\[
 G_\ell(m)=m+\sum_{j=1}^{\ell}\frac{m}{m+j},\qquad
 P_\ell(m)=m+\frac{m\ell}{m+\ell}\quad(m\ge1).
\]
We also set $G_\ell(0)=0$ for an empty null block.
The one-tail argument in the proof of \MainResult{thm:truncated}
gives, for arbitrary helpers,
\begin{equation}\label{supp:rankmanual-oracle}
 \int\psi(x)\,dx\le G_\ell(m).
\end{equation}

\begin{lemma}
\label{supp:rankmanual-onehelper}
Suppose that all but at most one of the $\ell\ge1$ false-null
p-values are nondecreasing functions of the tail distance $x$.
The remaining helper may be an arbitrary Borel function.  Then
\begin{equation}\label{supp:rankmanual-onehelper-bound}
 \FDR_{\mathrm{tail}}/h\le P_\ell(m).
\end{equation}
\end{lemma}

\begin{proof}
In the true-null-zeroed model of Lemma~\ref{supp:rankmanual-rank},
force the exceptional helper to equal one or zero and write the
resulting ranks as $r_-(x)$ and $r_+(x)$, respectively.  Both are
nonincreasing, and
\[
 m\le r_-(x)\le r_+(x)\le m+\ell.
\]
The actual effective rank equals one of these two ranks, rather than
merely lying between them. To verify this, let its value be $r$.
If the helper is rejected, its p-value is already at most $rh$;
lowering it to zero therefore leaves unchanged the counts at every
rank $k>r$, and preserves the crossing at $r$.  Hence $r=r_+$.
If it is not rejected, raising it to one leaves the count at $r$
unchanged and cannot create a larger crossing.  Hence $r=r_-$.
Equality at a critical value belongs to the rejected case.

Let $T_-$ and $T_+$ be the endpoints of the success intervals
of these two nonincreasing branches.  The endpoint argument in
Lemma~\ref{supp:rankmanual-rank} yields
\[
 m\le T_-\le T_+\le m+\ell,
 \qquad r_\pm(x)\ge T_\pm\quad(x<T_\pm).
\]
On $(0,T_-)$ either branch can reject the nulls, and the largest
possible FDP is $m/r_-$.  On $(T_-,T_+)$ only the $r_+$ branch
can reject them. After $T_+$ neither branch can reject them. Even allowing the
branch to be chosen freely at each point, we obtain
\[
 \begin{aligned}
 \int\psi(x)\,dx
 &\le\int_0^{T_-}\frac{m}{r_-(x)}\,dx
       +\int_{T_-}^{T_+}\frac{m}{r_+(x)}\,dx\\
 &\le m+m\frac{T_+-T_-}{T_+}
 \le m+m\left(1-\frac{m}{m+\ell}\right)=P_\ell(m).
 \end{aligned}
\]
Jumps and endpoint choices do not affect these integrals.  The
argument imposes no monotonicity on the exceptional helper.
\end{proof}

We next record the only cross-tail refinement needed in dimension
seven.  When $\ell=2$ and the helpers have opposite directions,
each tail has one nondecreasing helper and one other helper.  Thus
Lemma~\ref{supp:rankmanual-onehelper} gives
\begin{equation}\label{supp:rankmanual-pairbound}
 \FDR/h\le P_2(a)+P_2(b).
\end{equation}
The function $P_2$ is strictly concave, so for fixed $a+b$ its
integer maximum occurs at the most balanced split.

\begin{lemma}
\label{supp:rankmanual-asymmetric}
Suppose $\ell=2$, $1\le a<b$, and the two helpers have opposite
directions.  Then
\begin{equation}\label{supp:rankmanual-asymmetric-bound}
 \frac{\FDR}{h}\le
 \max\left\{a+G_2(b),\,
       P_2(a)+b\left(1+\frac1{b+2}\right)\right\}.
\end{equation}
\end{lemma}

\begin{proof}
Write $B(U)$ for the nondecreasing helper and $D(U)$ for the
nonincreasing one.  If $D(U)>(a+2)h$ throughout
$0<U\le(a+2)h$, it is unavailable at every possible left-tail
null rejection.  That tail contains only the nondecreasing helper
$B$ and contributes at most $ah$ by Lemma~\ref{supp:rankmanual-rank}.
The unrestricted right-tail bound \eqref{supp:rankmanual-oracle}
then gives the first term in \eqref{supp:rankmanual-asymmetric-bound}.

Otherwise there is a $u_0\le(a+2)h$ in that interval such that
$D(u_0)\le(a+2)h$.  At every possible right-tail null rejection,
\[
 U\ge1-(b+2)h\ge1-q>q\ge u_0.
\]
Global monotonicity of $D$ and the integer inequality $a+1\le b$
therefore imply
\[
 D(U)\le(a+2)h\le(b+1)h.
\]
The right-tail effective rank is consequently always at least
$b+1$, and at most $b+2$.  Its contribution is bounded by
\[
 h\left\{(b+1)\frac{b}{b+1}+\frac{b}{b+2}\right\}
 =hb\left(1+\frac1{b+2}\right).
\]
Applying Lemma~\ref{supp:rankmanual-onehelper} to the left tail gives
the second term. The two cases include arbitrary jumps: if a
single endpoint triggers the second case, monotonicity still gives
the bound throughout the right tail.
\end{proof}

\begin{proof}[Proof of \MainResult{rankmanual-main}]
Put $n=a+b$.  If $n=0$, the FDR is zero.  If $\ell=0$, the two
disjoint null-rejection events have total probability
$(a+b)h=q$.  If only one true-null block is present, then
\[
 \FDR/h\le G_\ell(n)
 =n+\sum_{j=1}^{\ell}\frac{n}{n+j}\le n+\ell=K.
\]
For two nonempty blocks and $\ell\ge1$, if all helpers have the
same direction, one tail is entirely monotone and the other can
be bounded without restrictions.  After exchanging the two tails
if necessary,
\begin{equation}\label{supp:rankmanual-samedirection}
 \FDR/h\le a+G_\ell(b)<a+b+\ell=K.
\end{equation}
In particular, this covers every $\ell=1$ case.

It remains to consider $a,b\ge1$, $\ell\ge2$, and mixed helper
directions.  The unrestricted two-tail bound is
$\FDR/h\le G_\ell(a)+G_\ell(b)$.  By concavity of $G_\ell$, its integer maximum at fixed $n$ is attained at
$(a,b)=(\lfloor n/2\rfloor,\lceil n/2\rceil)$.
The following exact scalar values give bounds for all remaining cases.
An asterisk marks a bound exceeding the target, which is $1$ for
$K\le7$ and $21/20$ for $K=8$.
\begin{center}
\small
\begin{tabular}{ccl}
\toprule
$K$ & True-null counts $n$ & Two-tail upper bounds for $\FDR/q$\\
\midrule
$2,3$ & --- & No remaining cases\\
$4$ & $2$ & $11/12$\\
$5$ & $2,3$ & $5/6,\ 1$\\
$6$ & $2,3,4$ & $137/180,\ 113/120,\ (19/18)^*$\\
$7$ & $2,3,4,5$ & $7/10,\ 53/60,\ (107/105)^*,\ (451/420)^*$\\
$8$ & $2,3,4,5,6$ & $363/560,\ 929/1120,\ 39/40,\ (101/96)^*,\ (87/80)^*$\\
\bottomrule
\end{tabular}
\end{center}

For $\ell=2$, mixed directions mean one helper of each direction,
so \eqref{supp:rankmanual-pairbound} applies.  Concavity gives bounds
$6$ and $42/5$ for $\FDR/h$ when $n=4$ and $n=6$, respectively.
When $n=5$, reflect the tails to arrange $a<b$ and use
Lemma~\ref{supp:rankmanual-asymmetric}.  The two possible splits give
\[
 \begin{array}{c|cc}
 (a,b)&a+G_2(b)&P_2(a)+b(1+1/(b+2))\\ \hline
 (1,4)&97/15&19/3\\
 (2,3)&127/20&33/5
 \end{array}
\]
and hence a common bound $33/5$.

For $\ell=3$, mixed directions have counts one and two.  In one
tail only one helper is not nondecreasing; use
Lemma~\ref{supp:rankmanual-onehelper} there and the unrestricted bound
on the other tail.  Maximizing over ordered positive splits gives
\[
 \max_{a+b=4}\{P_3(a)+G_3(b)\}=\frac{203}{30},\qquad
 \max_{a+b=5}\{P_3(a)+G_3(b)\}=\frac{121}{15}.
\]
These are comparisons of three and four scalar fractions, respectively;
the maximizing splits are $(2,2)$ and $(3,2)$. Dividing the five
remaining bounds by their respective dimensions gives at most $1$
for $K\le7$ and at most $21/20$ for $K=8$, as recorded in the main
text. The bound $121/120$ for $(K,n,\ell)=(8,5,3)$ is used only
to compare with $21/20$; it does not establish noninflation for
that fixed true-null count.

The all-null configuration with all p-values equal to $U$ attains
$q$ for every $K$.  For $K=8$, Theorem~\ref{supp:thm:two-tail} with $a=3$
and one helper pair ($m=1$ in that theorem) attains
\[
 \frac{\FDR}{q}=\frac68+\frac{12}{8\cdot5}=\frac{21}{20}.
\]
The upper bounds are therefore sharp, proving the theorem.
\end{proof}

\begin{corollary}
\label{supp:rankmanual-onefalse}
For $K\ge2$ and $0<q<1/2$,
\[
 \frac{f^{(1)}_{K-1,K}(q)}q=1-\frac1{K^2}.
\]
\end{corollary}

\begin{proof}
Put $n=K-1$.  There is only one helper, so its direction makes
one tail monotone.  Allowing an empty block to contribute zero,
the argument for \eqref{supp:rankmanual-samedirection} gives, after
reflection if needed,
\[
 \FDR/h\le a+G_1(b)=n+\frac{b}{b+1}
 \le n+\frac{n}{n+1}.
\]
For equality, take all $n$ true p-values equal to $U$ and the
single nonincreasing helper equal to one on $U\le nh$ and zero
on $U>nh$.  BH rejects the $n$ nulls on $U\le nh$, and all
$K$ hypotheses on $nh<U\le Kh=q$.  Hence
\[
 \FDR=nh+h\frac nK=q\left(1-\frac1{K^2}\right).
\]
\end{proof}


\begin{thebibliography}{}

\bibitem[Benjamini and Hochberg(1995)]{BH1995}
Benjamini, Y. and Hochberg, Y. (1995).
Controlling the false discovery rate: A practical and powerful approach to multiple testing.
\emph{Journal of the Royal Statistical Society: Series B (Methodological)}, \textbf{57}(1), 289--300.

\bibitem[Benjamini and Yekutieli(2001)]{BY2001}
Benjamini, Y. and Yekutieli, D. (2001).
The control of the false discovery rate in multiple testing under dependency.
\emph{The Annals of Statistics}, \textbf{29}(4), 1165--1188.

\bibitem[Block et~al.(2008)Block, Savits, and Wang]{BlockSavitsWang2008}
Block, H. W., Savits, T. H. and Wang, J. (2008).
Negative dependence and the Simes inequality.
\emph{Journal of Statistical Planning and Inference}, \textbf{138}(12), 4107--4110.

\bibitem[Chen et~al.(2026)Chen, Wang, Zhang, and Wu]{ChenWangZhangWu2026}
Chen, X., Wang, F., Zhang, H. and Wu, Z. (2026).
On correlation robustness of supremum-based p-value combination tests.
\emph{Statistica Sinica}, forthcoming.

\bibitem[Chen et~al.(2023)Chen, Liu, Tan, and Wang]{CLTW23}
Chen, Y., Liu, P., Tan, K. S. and Wang, R. (2023).
Trade-off between validity and efficiency of merging p-values under arbitrary dependence.
\emph{Statistica Sinica}, \textbf{33}(2), 851--872.

\bibitem[Chi et~al.(2025)Chi, Ramdas, and Wang]{ChiRamdasWang2024}
Chi, Z., Ramdas, A. and Wang, R. (2025).
Multiple testing under negative dependence.
\emph{Bernoulli}, \textbf{31}(2), 1230--1255.

\bibitem[Dobriban(2026)]{Dobriban2026}
Dobriban, E. (2026).
The Benjamini--Hochberg procedure can fail to control the FDR for correlated two-sided Gaussian tests.
\emph{Proceedings of the National Academy of Sciences}, \textbf{123}(41), e2626267123.

\bibitem[Guo and Bhaskara Rao(2008)]{GuoRao2008}
Guo, W. and Bhaskara Rao, M. (2008).
On control of the false discovery rate under no assumption of dependency.
\emph{Journal of Statistical Planning and Inference}, \textbf{138}(10), 3176--3188.

\bibitem[Hommel(1983)]{Hommel1983}
Hommel, G. (1983).
Tests of the overall hypothesis for arbitrary dependence structures.
\emph{Biometrical Journal}, \textbf{25}(5), 423--430.

\bibitem[Lei(2026)]{Lei2026}
Lei, L. (2026).
How much can Gaussian dependence inflate the Benjamini--Hochberg procedure's FDR?
\emph{arXiv preprint} 2607.14812v4.

\bibitem[Roquain and Villers(2011)]{RoquainVillers2011}
Roquain, {\'E}. and Villers, F. (2011).
Exact calculations for false discovery proportion with application to least favorable configurations.
\emph{The Annals of Statistics}, \textbf{39}(1), 584--612.

\bibitem[Samuel-Cahn(1996)]{SamuelCahn1996}
Samuel-Cahn, E. (1996).
Is the Simes improved Bonferroni procedure conservative?
\emph{Biometrika}, \textbf{83}(4), 928--933.

\bibitem[Simes(1986)]{Simes1986}
Simes, R. J. (1986).
An improved Bonferroni procedure for multiple tests of significance.
\emph{Biometrika}, \textbf{73}(3), 751--754.

\bibitem[Su(2018)]{Su2018}
Su, W. J. (2018).
The FDR-linking theorem.
\emph{arXiv preprint} 1812.08965.

\bibitem[Vovk et~al.(2022)Vovk, Wang, and Wang]{VWW22}
Vovk, V., Wang, B. and Wang, R. (2022).
Admissible ways of merging p-values under arbitrary dependence.
\emph{The Annals of Statistics}, \textbf{50}(1), 351--375.

\bibitem[Zhang(2026)]{Zhang2026}
Zhang, X. (2026).
Exact Simes inequalities for two-sided Gaussian tests.
\emph{arXiv preprint} 2611.02513.

\end{thebibliography}
\end{document}